\documentclass[11pt, reqno]{amsart}
\usepackage{amsmath,amsthm,amsfonts,amssymb,mathrsfs,bm,graphicx,stmaryrd}
\usepackage{subcaption}\usepackage{mathtools}
\usepackage{dsfont}
\usepackage{multicol}
\usepackage[colorlinks=true,linkcolor=blue,citecolor=blue]{hyperref}
\usepackage{xcolor}
\usepackage{bbm}

\newcommand{\QQ}{\mathbb{Q}}
\newcommand{\wcL}{\widetilde{\mathcal{L}}}
\newcommand{\res}{\mathrm{res}}

\newcommand{\cI}{\mathcal{I}}
\newcommand{\RR}{\mathbb{R}}

\newcommand{\PP}{\mathbb{P}}
\newcommand{\EE}{\mathbb{E}}

\newcommand{\0}{\mathbf{0}}

\newcommand{\boxx}{\mathrm{Box}}

\newcommand{\diam}{\mathrm{diam}}

\newcommand{\wT}{\widetilde{T}}

\newcommand{\leb}{\mathrm{Leb}}

\newcommand{\wX}{\widetilde{X}}

\newcommand{\NN}{\mathbb{N}}
\newcommand{\NU}{\mathrm{NU}}
\usepackage[letterpaper,hmargin=1in,vmargin=1in]{geometry}
\parskip 	\smallskipamount

\newtheorem{theorem}{Theorem}
\newtheorem{lemma}[theorem]{Lemma}

\newtheorem{corollary}[theorem]{Corollary}
\newtheorem{proposition}[theorem]{Proposition}

\theoremstyle{definition}

\newcommand{\fF}{\mathfrak{F}}
\newcommand{\bv}{\mathbf{v}}
\newcommand{\bw}{\mathbf{w}}
\newcommand{\cA}{\mathcal{A}}
\newcommand{\cB}{\mathcal{B}}

\newcommand{\cC}{\mathcal{C}}

\newcommand{\cE}{\mathcal{E}}
\newcommand{\cN}{\mathcal{N}}

\newcommand{\ff}{\mathfrak{f}}
\newcommand{\cW}{\mathcal{W}}
\newcommand{\cL}{\mathcal{L}}

\newcommand{\three}{\mathfrak{T}}

\newcommand{\cS}{\mathcal{S}}

\newcommand{\inte}{\mathrm{int}}

\newcommand{\ZZ}{\mathbb{Z}}
\newcommand{\dis}{\mathrm{dis}}

\DeclareMathOperator*{\argmax}{\arg\! \max}

\newcommand{\cH}{\mathcal{H}}
\newcommand{\cT}{\mathcal{T}}
\usepackage[backend=biber,doi=false,style=alphabetic,url=false,maxalphanames=10,maxnames=50]{biblatex}
\AtBeginBibliography{\small}

\begin{document}
\title[]{Duality in the directed landscape and its applications to fractal geometry}
\author[]{Manan Bhatia}
\address{Manan Bhatia, Department of Mathematics, Massachusetts Institute of Technology, Cambridge, MA, USA}
\email{mananb@mit.edu}
\date{}
\maketitle
\begin{abstract}
 Geodesic coalescence, or the tendency of geodesics to merge together, is a hallmark phenomenon observed in a variety of planar random geometries involving a random distortion of the Euclidean metric. As a result of this, the union of interiors of all geodesics going to a fixed point tends to form a tree-like structure which is supported on a vanishing fraction of the space. Such geodesic trees exhibit intricate fractal behaviour; for instance, while almost every point in the space has only one geodesic going to the fixed point, there exist atypical points which admit two such geodesics. In this paper, we consider the directed landscape, the recently constructed \cite{DOV18} scaling limit of exponential last passage percolation (LPP), with the aim of developing tools to analyse the fractal aspects of the tree of semi-infinite geodesics in a given direction. We use the duality \cite{Pim16} %
 between the geodesic tree and the interleaving competition interfaces in exponential LPP to obtain a duality between the geodesic tree and the corresponding dual tree in the landscape. %
Using this, we show that problems concerning the fractal behaviour of sets of atypical points for the geodesic tree can be transformed into corresponding problems for the dual tree, which might turn out to be easier. In particular, we use this method to show that the set of points admitting two semi-infinite geodesics in a fixed direction a.s.\ has Hausdorff dimension $4/3$, thereby answering a question posed in \cite{BSS22}. We also show that the set of points admitting three semi-infinite geodesics in a fixed direction is a.s.\ countable.
\end{abstract}

\begin{figure}
  \centering
  \includegraphics[width=0.49\textwidth]{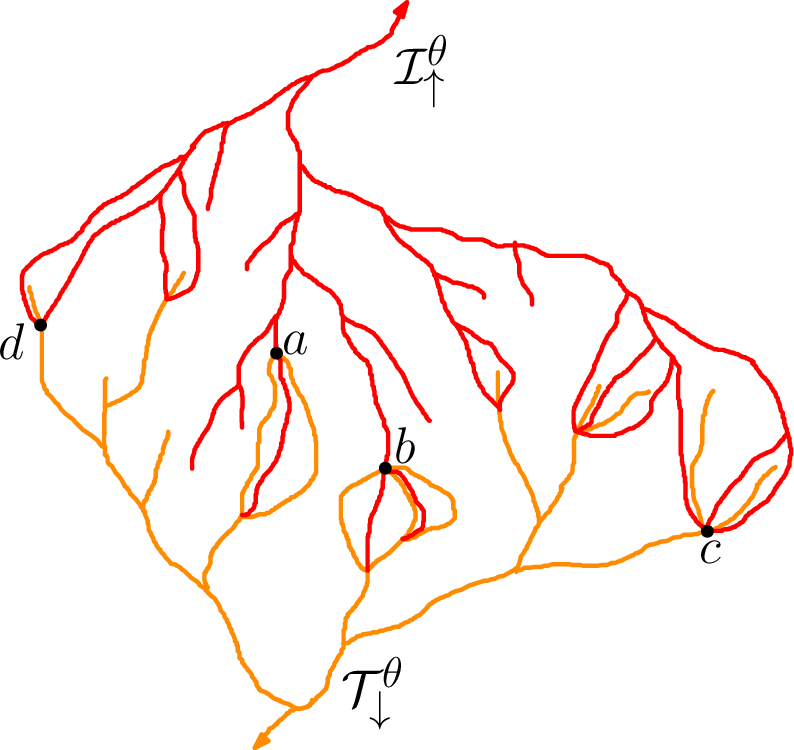}
  \caption{The downward $\theta$-directed geodesic tree $\cT^\theta_\downarrow$ and its corresponding interface portrait (dual tree) $\cI^\theta_\uparrow$. As can be seen by looking at the points $a$ and $d$, points on one tree equal the non-uniqueness points corresponding to the other tree. Similarly, as witnessed by the points $b$ and $c$, the trifurcation points of one tree are exactly the points with three distinct paths emanating out in the other tree. Note that $a,b\notin \cT_\downarrow^\theta$ and $c,d\notin \cI_\uparrow^\theta$ since ``leaves'' are not considered a part of the tree.}
  \label{fig:cover}
\end{figure}

\section{Introduction}
\label{sec:intro}
In the past few years, random geometry has emerged as a subject of its own following major developments in the study of certain natural continuum models-- the directed landscape \cite{DOV18}, the Liouville quantum gravity metrics \cite{DDDF20,GM21} along with the specially integrable case of Brownian geometry \cite{LeGal10}. While these models have very different origins, they are rather similar in the sense that they comprise of a random distance function whose geodesics tend to coalesce with each other and display rich fractal behaviour, and the study of this fractal behaviour has recently attracted significant interest in all three cases-- the directed landscape (e.g.\ \cite{GZ22, CHHM21}), Brownian geometry (e.g.\ \cite{Mie13, LeGal22, MQ20}) and Liouville quantum gravity (e.g.\ \cite{Gwy21,FG22}). %

Instead of looking at just one geodesic at once, one can consider the set of all geodesics going to a fixed point and attempt to study the fractal geometry of this even richer object, which, by geodesic coalescence, often turns out to be a tree. Indeed, the primary tool \cite{CS04, LeGal07} for studying the Brownian map is to consider such a  geodesic tree to a typical point and obtain a strong notion of integrability for its ``dual'', by which we mean the set of continuous curves ``avoiding'' the geodesic tree. This dual object turns out to be the celebrated continuum random tree \cite{Ald91} of Aldous along with an independent Brownian motion on it which encodes the distances of general points to the aforementioned typical point.
Drawing inspiration from the above, the aim of this work is to develop a ``duality'' for the directed landscape which gives us the ability to convert an atypical set of points for the geodesic tree into a corresponding set of points for the ``dual'' forest, which we show is a tree and call as the interface portrait, extending terminology used in \cite{RV21}. As we demonstrate, some sets of atypical points turn out to become much more tractable in the dual picture, thereby allowing us to understand their fractal structure. %

Since the route to obtaining duality in the continuum will be via a corresponding duality present in the discrete \cite{Pim16,FMP09,SS21+}, we begin by introducing exponential last passage percolation (LPP), an integrable planar lattice model which converges \cite{DV21} to the directed landscape in the scaling limit. Consider the lattice $\ZZ^2$ and endow each vertex $u\in \ZZ^2$ with an i.i.d.\ $\exp(1)$ weight $X_u$. For any two points $u,v$ which are ordered in the sense that $v\geq u$ coordinate-wise, we define the passage time $T(u,v)$ to be the largest weight of an up-right lattice path from $u$ to $v$, where the weight $\ell(\pi)$ of such a lattice path $\pi$ is defined to be $\ell(\pi)=\sum_{w\in \pi\setminus \{u\}}X_w$. Interpreting $T(u,v)$ as ``distances''\footnote{We caution that due to the maximization present in the definition of path lengths, $T(u,v)$ actually yields an anti-metric in the sense that it satisfies the triangle inequality in reverse. The same is true for the directed landscape.} in a directed space, we can think of exponential LPP as a discrete model of random geometry.

Having defined LPP, we now introduce the directed landscape, a continuum model of random geometry introduced in \cite{DOV18} and known \cite{DV21} to be the space-time scaling limit of exponential LPP. Consider the space
\begin{equation}
  \label{eq:36}
  \RR_\uparrow^4=\{(x,s;y,t)\in \RR^4:s<t\}.
\end{equation}
The directed landscape $\cL$ is a random real valued function on $\RR_\uparrow^4$ satisfying, for any $s<r<t$, the natural composition law
\begin{equation}
  \label{eq:9}
  \cL(x,s;y,t)=\max_{z\in \RR} \{ \cL(x,s;z,r) + \cL( z,r;y,t)\},
\end{equation}
with the value $\cL(x,s;y,t)$ to be interpreted as the last passage time from the point $(x,s)$ to $(y,t)$. 
With the notion of passage times at hand, one can further define the weight $\ell(\eta)$ of a continuous function $\eta\colon [s,t]\rightarrow \RR$, thereafter referred to as a path, by
\begin{equation}
  \label{eq:37}
  \ell(\eta)= \inf_{k\in \NN} \inf_{s=t_0<t_1<\cdots<t_k=t}\sum_{i=1}^k\cL(\eta(t_{i-1}),t_{i-1};\eta(t_i),t_i).
\end{equation}
We will always think of paths $\eta\colon [s,t]\rightarrow \RR$ as drawn vertically in the plane $\RR^2$ with the vertical coordinate being thought of as the time direction. In line with this, the graph of a path $\eta$ will denote the set $\{(\eta(r),r): r\in [s,t]\}\subseteq \RR^2$.

One can now consider all such paths $\eta\colon [s,t]\rightarrow \RR$ with $\eta(s)=x$ and $\eta(t)=y$ and maximise $\ell(\eta)$ over such paths, and any path achieving the above-mentioned maximum is called a geodesic and is denoted as $\gamma_{(x,s)}^{(y,t)}$. As one would expect, geodesics a.s.\ exist \cite{DOV18} between any two ordered points and for fixed points $(x,s),(y,t)$ with $s<t$, there is almost surely a unique geodesic $\gamma_{(x,s)}^{(y,t)}$. Apart from finite geodesics, it can in fact be shown that almost surely, simultaneously for every direction $\theta\in \RR$ and every point $(y,t)\in \RR^2$, there exists \cite{BSS22} (see also \cite{RV21,GZ22}) a semi-infinite geodesic $\Gamma_{(y,t)}^\theta\colon (-\infty,t]\rightarrow \RR$ satisfying $\Gamma_{(y,t)}^\theta(s)/s\rightarrow \theta$ as $s\rightarrow -\infty$, in the sense that any finite segment of $\Gamma_{(y,t)}^\theta$ is a geodesic. Again, for any fixed $\theta\in \RR$ and $(y,t)\in \RR^2$, one a.s.\ has a unique semi-infinite geodesic $\Gamma_{(y,t)}^\theta$, but there might exist atypical points $(y,t)$ for which this fails. In fact, for a fixed $\theta\in \RR$, the union of the graphs of $\Gamma^\theta_{(y,t)}\lvert_{(-\infty,t)}$ over all such semi-infinite geodesics and all points $(y,t)\in \RR^2$ is usually called (e.g.\ \cite[Section 3.7]{RV21}) the geodesic tree $\cT_\downarrow^\theta$, since the corresponding object is easily seen to be a tree in the prelimit; we give a formal definition of the set $\cT_\downarrow^\theta$ in \eqref{eq:42}. It is indeed true that, just as in the discrete, the above continuum object $\cT_\downarrow^\theta$ also forms a tree as we shall show in Corollary \ref{thm:3}, thereby answering a question raised in \cite{BSS22}.

The geometry of the geodesic tree $\cT^\theta_\downarrow$ displays intricate fractal behaviour, of which many aspects are not well understood. Investigating such fractal behaviour follows the general theme of studying fractality in random geometry, and this is a topic which has received significant interest recently for the directed landscape. Some recent successes include the study of atypical stars \cite{BGH21,BGH19,GZ22,Bha22}, the analysis of exponents governing $k$ disjoint geodesics starting and ending nearby \cite{Ham20}, the study of exceptional times where the Kardar-Parisi-Zhang (KPZ) fixed point has multiple maxima \cite{CHHM21,Dau22}, and the description of the atypical set of directions admitting non-coalescing semi-infinite geodesics \cite{BSS22}. We refer the reader to \cite{GH22} for a survey on fractality in the directed landscape.

The geodesic tree $\cT^\theta_\downarrow$ leads to many natural collections of corresponding atypical points, and in this paper, we attempt to analyse two such sets-- the set of points having two semi-infinite geodesics in a fixed direction $\theta$, and the set of points admitting three such geodesics. The question of computing the fractal dimension of the former was raised in \cite{BSS22}, and in fact it was not even known if points of the latter type exist. In this work, we successfully solve the above two questions, showing that the former set a.s.\ has Hausdorff dimension $4/3$, while the latter is a.s.\ countable, where by `countable', we mean countably infinite. To do so, we consider the duality \cite{Pim16,FMP09} between the geodesic tree and the corresponding competition interfaces in exponential LPP and show that it passes unscathed to the scaling limit. Using the convergence properties of the directed landscape \cite{DV21}, this manifests in there being a notion of a dual landscape $\wcL$ coupled to a directed landscape $\cL$ such that the geodesic tree and interface portrait interchange when one passes between the primal and dual picture. This of course implies that the interface portrait and the geodesic tree have the same distribution up to a reflection and this allows us to track how atypical points for the geodesic tree appear in the interface portrait picture, with the hope that the dual picture will be more amenable to understanding certain fractal aspects. We hope that the above technique of tracking the behaviour of atypical points on passage to the dual picture will be helpful in investigating further fractal properties of the geodesic tree in the future.

In fact, as described to us by B{\'a}lint Vir{\'a}g and an anonymous referee after the first version of this paper was posted online, the setting and techniques of this paper have a striking parallel in the well-studied object known as the `Brownian Web'. First conceived in \cite{Arr79, Arr81} and recently surveyed in \cite{SSS17}, the Brownian web, or the Arratia flow, consists of the tree formed by a system of independent coalescing one dimension Brownian motions started from every point in $\RR^2$, with the vertical axis viewed as the time axis. Just as the interface portraits $\cI_\downarrow^\theta$ can be defined as the dual objects to the geodesic trees $\cT_\uparrow^\theta$ (see \eqref{eq:15}), one can similarly consider dual \cite{Arr81,TW98,FINR06} semi-infinite paths for the Brownian web, and these are known to form a tree as well. Further, up to a reflection, this dual tree has the same law as the Brownian web, an aspect perfectly in line with the setting of this paper. Finally, similar to Theorem \ref{thm:6}, the above-mentioned duality along with the interlocked nature of the primal and dual Brownian webs have been used to investigate (\cite[Proposition 2.4]{TW98}, \cite[Theorems 3.11-3.14]{FINR06}) the fractal aspects of atypical points for the Brownian web. As a result of the above similarities, the interlocked primal and dual geodesic trees exhibited in this paper can together be viewed as a KPZ analogue of the Brownian web, where the individual semi-infinite paths now display the characteristic KPZ $2/3$ wandering exponent, as opposed to the $1/2$ exponent observed in the case of Brownian motion.
\section{Main results}
\label{sec:results}
\subsection{Geodesics do not form bubbles}
\label{sec:later}
The following independently interesting structural result about geodesics is easy to state but has significant consequences; we note that a similar result recently found application \cite{MQ20} in the allied field of Brownian geometry in the proof of a strong form of coalescence of geodesics. %

\begin{theorem}
  \label{thm:2}
  Almost surely, there does not exist any $u=(x,s;y,t)\in \RR_\uparrow^4$ such that there are two distinct geodesics $\eta^1,\eta^2$ from $(x,s)$ to $(y,t)$ with the property that for some small enough $\delta>0$, $\eta^1(r)=\eta^2(r)$ for all $r\in(s,s+\delta)\cup (t-\delta,t)$. As a consequence, almost surely, for any geodesic $\gamma\colon[a,b]\rightarrow\RR$ and any $(a_1,b_1)\subseteq [a,b]$, $\gamma\lvert_{(a_1,b_1)}$ is the unique geodesic between its endpoints.
\end{theorem}
The configurations which the above result denies the existence of, will be referred to as geodesic bubbles, and we point the reader to Figure \ref{fig:bubble} for an illustration. The proof of Theorem \ref{thm:2} is delicate since it concerns all possible geodesics, and these might a priori look very different from geodesics between typical points. %
The proof proceeds by first using a rational approximation argument along with a soft application of the recent result \cite{Bha22} concerning the fractality of atypical points on geodesics (see Proposition \ref{prop:8}), to obtain that there must be many ``typical'' points on any geodesic. After this, one can use another rational approximation argument to show that the existence of a bubble as in Theorem \ref{thm:2} would contradict the uniqueness of geodesics between rational points. We refer the reader to Figure \ref{fig:bubble} for a summary of the latter part of the above proof.

In order to discuss some consequences of the above result, we now introduce semi-infinite geodesics, and as we shall see soon, these will play a much bigger role in this work than finite geodesics. For any $p=(x,s)\in \RR^2$, an upward semi-infinite geodesic emanating from $p$ is a path $\gamma \colon [s,\infty)\rightarrow \RR$ with the property that $\gamma(s)=x$ and for each $t>s$, $\gamma\lvert_{[s,t]}$ is a finite geodesic. If the path $\gamma$ above satisfies $\lim_{t\rightarrow \infty} \gamma(t)/t=\theta$ for some $\theta\in \RR$, then we call it an upward $\theta$-directed semi-infinite geodesic emanating from $p$ and use the notation $\Gamma_{p,\uparrow}^\theta$ to refer to it. We note that one can analogously define semi-infinite geodesics and $\theta$-directed semi-infinite geodesics in the downward direction. We denote the latter by $\Gamma^\theta_{p,\downarrow}$ and these satisfy $\lim_{t\rightarrow -\infty}\Gamma_{p,\downarrow}^\theta(t)/t=\theta$. We will often abbreviate $\Gamma^\theta_{p,\downarrow}$ to just $\Gamma^\theta_p$ and further use $\Gamma_p$ for $\Gamma^0_p$. As we shall see soon, we will frequently work with the interior $\inte (\Gamma_p^\theta)$ of a geodesic $\Gamma_p^\theta$, where by the interior of a path or a semi-infinite path, we simply mean its graph with its endpoints removed.

Before continuing the discussion on semi-infinite geodesics, we make a notational remark. Just as we have downward and upward semi-infinite geodesics in the above paragraph, we will have both downward and upward versions of many objects throughout the paper. Owing to the flip symmetry of the directed landscape (see Proposition \ref{prop:symm}), we will often state results only for one of these versions, and it will be understood that the same result applies to the other version as well.

In \cite[Theorem 2.5 (i)]{BSS22}, it was shown that almost surely, every downward semi-infinite geodesic must be $\theta$-directed for some $\theta \in \RR$, and further, such $\theta$-directed downward semi-infinite geodesics $\Gamma_{p}^\theta$ were shown to exist simultaneously for each direction $\theta$ and point $p$ (for semi-infinite geodesics in LPP, see \cite{FMP09,FP05}), with there always \cite[Theorem 6.5 (i)]{BSS22} being a uniquely defined left-most such geodesic $\underline{\Gamma}_p^\theta$ and right-most such geodesic $\overline{\Gamma}_p^\theta$. Further, it is known \cite{RV21,GZ22} that a.s.\ for each fixed direction $\theta$, $\Gamma^\theta_p$ must all coalesce in the sense that for any two geodesics $\Gamma^\theta_{p},\Gamma^\theta_{q}$, we have $\Gamma^\theta_{p}(t)=\Gamma^\theta_{q}(t)$ for all negative enough $t$. In fact, in \cite{BSS22}, the above was extended to hold for all $\theta\in \Xi_\downarrow^c$, where $\Xi_\downarrow$ is a random countable set of exceptional directions, which we call the directions of non-uniqueness. %
However, as described in \cite{BSS22}, the above coalescence along the uniqueness directions $\Xi_\downarrow^c$, is still not sufficient to establish that for all $\theta\in \Xi_\downarrow^c$, the set
\begin{equation}
  \label{eq:42}
\cT^\theta_\downarrow=\bigcup_{p\in\RR^2}\inte (\Gamma^\theta_p),
\end{equation}
where the union is taken over all choices of geodesics $\Gamma^\theta_p$, forms a tree, where throughout the paper, we call a subset $\mathrm{T}\subseteq \RR^2$ a tree if for any two points $p,q\in \mathrm{T}$, there is a unique (up to monotone reparametrization) $\mathrm{T}$-valued injective continuous curve from $p$ to $q$. Note the important fact that $\cT_\downarrow^\theta\neq \RR^2$ a.s.\ for all $\theta \in\Xi_\downarrow^c$, and this is shown later in Lemma \ref{lem:1}. Thus, it is crucial to use $\inte(\Gamma_p^\theta)$ in \eqref{eq:42} and indeed, $\bigcup_{p\in \RR^2}\Gamma_p^\theta$ is a.s.\ simply equal to $\RR^2$. We note that a definition similar to \eqref{eq:42} can be made for $\cT_\uparrow^\theta$ by considering upward geodesics. 

To understand what might possibly prevent $\cT_\downarrow^\theta$ from forming a tree, we consider the following sets defined in \cite{BSS22},
\begin{align}
  \label{eq:10}
  \NU_0^\theta&=
  \left\{
  (y,t)\in \RR^2: \underline{\Gamma}_{(y,t)}^\theta(s)\neq \overline{\Gamma}_{(y,t)}^\theta(s)\textrm{ for some }s<t
  \right\},\\
  \NU_1^\theta&=
        \left\{
        (y,t)\in \RR^2: \exists ~\varepsilon>0 \textrm{ such that } \underline{\Gamma}_{(y,t)}^\theta(s)\neq \overline{\Gamma}_{(y,t)}^\theta(s)\textrm{ for all } s\in (t-\varepsilon,t)
        \right\}.
\end{align}
That is, while points in both sets exhibit non-unique downward $\theta$-directed semi-infinite geodesics, the latter set consists of points where the left-most and right-most such geodesics split immediately, while in the former case, they might possibly stay together for some time before splitting. We note that though $\NU_0^{\theta},\NU_1^\theta$ concern downward geodesics, we can similarly define the corresponding sets for upward geodesics and we denote these by $\NU_{0,\uparrow}^\theta,\NU_{1,\uparrow}^\theta$. Though it is clear that $\NU^\theta_1\subseteq \NU^\theta_0$ for all $\theta\in \Xi_\downarrow^c$, it is a priori possible that the inclusion is strict. Note that the existence of $p\in \NU^\theta_0\setminus \NU_1^\theta$ for some $\theta\in \Xi_\downarrow^c$ would create a cycle in $\cT_\downarrow^\theta$ and imply that $\cT_\downarrow^\theta$ cannot be a tree. In fact, the validity of $\NU^\theta_0=\NU^\theta_1$ for the uniqueness directions $\theta\in \Xi_\downarrow^c$ was listed as an open problem in \cite{BSS22}, and we obtain it as a corollary of Theorem \ref{thm:2}. 
\begin{corollary}
  \label{thm:3}
  Almost surely, simultaneously for all directions $\theta\in \Xi_\downarrow^c$, we have the equality $\NU^\theta_0=\NU^\theta_1$. Consequently, almost surely, simultaneously for every $\theta\in \Xi_\downarrow^c$, $\cT_\downarrow^\theta$ is a one-ended tree in the sense that it is a tree in which any two semi-infinite paths eventually merge and stay together.
\end{corollary}
Indeed, if the inclusion $\NU_{1}^\theta\subseteq \NU_0^\theta$ were strict, then we would have the existence of geodesic bubbles which is precluded by Theorem \ref{thm:2}. Before moving on, we mention that Theorem \ref{thm:2} has an interesting consequence (Proposition \ref{lem:50}) which states that the union of interiors of all geodesics is equal to the union of interiors of geodesics between rational points. We note that this upgrades the result \cite[Corollary 3.5]{DSV22} which states the same but only for left-most and right-most geodesics instead of all geodesics.
\subsection{The duality between geodesic trees and interface portraits}
\label{sec:duality}

In last passage percolation, competition interfaces \cite{FMP06} arise as the regions of space lying on the interface between two competing growth clusters. %
In fact, competition interfaces can also be interpreted \cite{Pim16} as semi-infinite paths living on the dual lattice which lie in between semi-infinite geodesics. The connection of competition interfaces to second class particles \cite{FP05} in the totally asymmetric simple exclusion process (TASEP) and their distributional equality \cite{FMP09} to semi-infinite geodesics has long been studied.

Extending the notion of a competition interface to the continuum limit, an object called the interface portrait, consisting of interfaces (or, the dual paths) between geodesics, was studied in \cite{RV21} (see also \cite{GZ22}) for the directed landscape coupled with an initial condition. Given a nice enough initial condition $f$ on $\RR\times \{0\}$ along with the directed landscape restricted to $\RR\times [0,\infty)$, \cite{RV21} defined a collection of interface paths emanating from $\RR\times \{0\}$ which lie ``between'' the set consisting of the interiors of geodesics from points in $\RR\times (0,\infty)$ to the initial condition $f$. We will use some results from the above work in this paper, and these will be discussed in Section \ref{sec:intface}.

In this paper, we extend the above-mentioned usage of the term `interface portrait' to refer to the set $\cI_\downarrow^\theta$ (resp.\ $\cI_\uparrow^\theta$) formed by the dual paths which interlace between the paths of the geodesic tree $\cT_\uparrow^\theta$ (resp.\ $\cT_\downarrow^\theta$). %
For $\theta\in \Xi_\uparrow^c$, we now define the $\theta$-directed interface portrait $\cI_\downarrow^\theta$ by
\begin{equation}
  \label{eq:15}
  \cI_\downarrow^\theta\coloneqq \bigcup_{\pi: \inte(\pi)\cap \cT^\theta_\uparrow=\emptyset}\inte(\pi),
\end{equation}
where the above union is over all downward semi-infinite paths $\pi$ satisfying $\inte(\pi)\cap \cT^\theta_\uparrow=\emptyset$. Note the obvious yet important property $\cI_\downarrow^\theta \cap \cT_\uparrow^\theta=\emptyset$. %
An advantage of the definition \eqref{eq:15} is that it canonically associates the interface portraits $\cI_\downarrow^\theta$ to the directed landscape a.s.\ simultaneously for all directions $\theta\in \Xi_\uparrow^c$ as opposed to a fixed direction $\theta\in \RR$. We note that as opposed to the setting in \cite{RV21}, the above-mentioned interface portraits are associated just to the directed landscape $\cL$, instead of a directed landscape coupled with an initial condition.

Often, it will be useful to consider the individual semi-infinite paths appearing in the union \eqref{eq:15}, and we now introduce notation for this. For any point $p=(x,s)$, if we have a semi-infinite path $\pi\colon (-\infty,s]\rightarrow \RR$ with $\pi(s)=x$ and $\inte(\pi)\cap \cT_\uparrow ^\theta=\emptyset$, then we call $\pi$ a downward $\theta$-directed interface emanating from $p$ and denote it by $\Upsilon_p^\theta$. Later, in Section \ref{sec:intface}, we shall see that almost surely, interfaces $\Upsilon_p^\theta$ exist for all points $p\in \RR^2$ and all directions $\theta \in \Xi_\uparrow^c$, and further, more than one choice of the aforementioned interface exists if and only if (see Lemma \ref{lem:53}) $p\in \cT_\uparrow^\theta$, and in particular (see Lemma \ref{lem:1}) for a fixed point $p$, there is a.s.\ a unique interface $\Upsilon_p^\theta$.

Now, with the above definition of the interfaces $\Upsilon_p^\theta$ at hand, we have the equality
\begin{equation}
  \label{eq:23}
  \cI_\downarrow^{\theta}=\bigcup_{p\in \RR^2}\inte(\Upsilon_{p}^{\theta}),
\end{equation}
where the union is over all possible interfaces in the case of non-uniqueness. %
That is, $\cI^\theta_\downarrow$ is the union of the interiors of all downward $\theta$-directed interfaces. Though we do not expand upon this, we can analogously define upward interfaces $\Upsilon^{\theta}_{p,\uparrow}$ and the corresponding interface portrait $\cI^\theta_\uparrow$ by just replacing the upward geodesic tree $\cT^{\theta}_{\uparrow}$ with the downward geodesic tree $\cT^{\theta}_\downarrow$ in the definitions. This will be important in the upcoming Section \ref{sec:duality}. To summarize, the interface portrait $\cI_\downarrow^\theta$, which will be shown to be a tree later, consists of dual paths (see Figure \ref{fig:cover}) interlacing between the gaps of the geodesic tree $\cT^\theta_\uparrow$. Finally, we note that we will sometimes work together with multiple directed landscapes and in such situations, we will use notation of the form $\cT^\theta_\downarrow(\cL)$ (or $\cI^\theta_\downarrow(\cL)$) to make explicit that the geodesic tree (or interface portrait) is being considered with respect to the landscape $\cL$.

For exponential LPP, the work \cite{Pim16} proved an exact duality for the geodesic tree and the corresponding dual tree consisting of competition interfaces. Precisely, it was shown that in exponential LPP with the weights $X=(X_{i,j})_{(i,j)\in \ZZ^2}$, the tree $\cT_\uparrow(X)$ consisting of semi-infinite geodesics going in the direction $(1,1)$ has the same distribution as its dual tree $\cI_\downarrow(X)\subseteq (\ZZ^2)^*=\ZZ^2+(1/2,1/2)$ up to a reflection, where we note that we have overloaded the notation $\cI,\cT$ by using it for LPP as well. Moreover (see Proposition \ref{prop:7}), it was established \cite{Pim16, FMP09} that for each choice of the LPP weights $X=\{X_{(i,j)}\}_{(i,j)\in \ZZ^2}$, there is a dual LPP model $\wX=\{\wX_{(i^*,j^*)}\}_{(i^*,j^*)\in (\ZZ^2)^*}$ with the property that $\cT_\downarrow(X)=\cI_\downarrow(\wX)$ and $\cI_\uparrow(X)=\cT_\uparrow(\wX)$.
In other words, the duality operation  transforms the noise $X$ into a dual noise $\widetilde{X}$ in a manner which swaps interface portraits and geodesic trees.

Though the equality $\Upsilon_p\stackrel{d}{=}\Gamma_p$ in the directed landscape was established in \cite{RV21,GZ22} for any fixed point $p$, the above striking results in the prelimit beg the question of whether we in fact have a similar distributional equality for the trees $\cT_\downarrow$ and $\cI_\downarrow$ as well. %
The following result establishes the above along with a version of the prelimiting duality $X\leftrightarrow \wX$ for the directed landscape. 

\begin{theorem}
  \label{thm:20}
  For any fixed direction $\theta\in \RR$, there is a coupling $(\cL_\theta,\wcL_\theta)$ of directed landscapes for which we have the almost sure equalities $\cT^\theta_{\downarrow}(\cL_\theta)=\cI^\theta_{\downarrow}(\wcL_\theta)$ and $\cI^\theta_{\uparrow}(\cL_\theta)=\cT^\theta_{\uparrow}(\wcL_\theta)$. %
\end{theorem}
In our construction of the above coupling, we do expect that $\cL_\theta$ and $\wcL_\theta$ are measurable with respect to each other just as $X$ and $\wX$ determine each other. However, we are currently not able to establish this. %

\subsection{The interface portrait $\cI_\downarrow ^\theta$ is a one-ended tree simultaneously for all directions $\theta\in \Xi_\uparrow^c$}
\label{sec:bigeod}

As we shall see later, it can be shown by the results of \cite{RV21} that for any $\theta\in \Xi_\uparrow^c$, and any $p,q\in \RR^2$, any two interfaces $\Upsilon_{p}^{\theta},\Upsilon_{q}^{\theta}$ must coalesce if they meet, and it is also not difficult to show that $\cI^\theta_\downarrow$ does not have any cycles (see Lemma \ref{lem:56}), and later, we shall refer to the above property as $\cI^\theta_\downarrow$ being a downward directed forest simultaneously for all $\theta\in \Xi_\uparrow^c$.%

 We note that the duality in Theorem \ref{thm:20} immediately improves the above to the statement that $\cI_\downarrow^\theta$ is a.s.\ a one-ended tree for any fixed direction $\theta$, since the geodesic tree $\cT^\theta_\downarrow$ does satisfy this almost surely. However, we can in fact go further, and show that $\cI_\downarrow ^\theta$ is not just a forest, but actually a one-ended tree a.s.\ simultaneously for all directions $\theta\in \Xi_\uparrow^c$.

\begin{theorem}
  \label{thm:4}
  Almost surely, the interface portraits $\cI_\downarrow^\theta$ are one-ended trees simultaneously for all directions $\theta\in \Xi_\uparrow^c$.
\end{theorem}

To show Theorem \ref{thm:4} (see Figure \ref{fig:int-tree}), we first show (Proposition \ref{prop:5}) that the directed landscape cannot have a bi-infinite geodesic, that is, there does not exist any bi-infinite path $\Gamma\colon (-\infty,\infty)\rightarrow \RR$ such that $\Gamma\lvert_{[s,t]}$ is a geodesic for every $s<t$. Such a statement has been shown \cite{BHS22,Pim16,BSS20} for the prelimiting model of exponential LPP, and we adapt the argument from \cite{BHS22} to the setting of the directed landscape. To go from the non-existence of bi-infinite geodesics to Theorem \ref{thm:4}, we use a compactness argument relying on the fact that $\cT^\theta_\uparrow$ and $\cI^\theta_\downarrow$ cannot intersect each other.

In view of the one-endedness from Theorem \ref{thm:4}, it will be useful in the next section to have the following definition of trifurcation points. For $\theta\in \Xi_\downarrow^c$, we say that a point $v\in \mathrm{Tri}(\cI^\theta_\uparrow)\subseteq \RR^2$ if and only if there exist some interfaces $\Upsilon^\theta_{p,\uparrow},\Upsilon^\theta_{q,\uparrow}$, starting from some points $p,q\in \RR^2$, which meet for the first time at $v$. Similarly, for $\theta\in \Xi_\uparrow^c$, we define $\mathrm{Tri}(\cT^\theta_\uparrow)$ such that $v\in \mathrm{Tri}(\cT^\theta_\uparrow)$ if and only if there exist some geodesics $\Gamma^\theta_{p,\uparrow},\Gamma^\theta_{q,\uparrow}$, emanating from some points $p,q\in \RR^2$, which meet for the first time at $v$.

\subsection{Investigating the fractal behaviour of atypical points via duality}
\label{sec:dual-except}
The duality in Theorem \ref{thm:20} provides a strong tool to study the fractal aspects of the geodesic tree in a fixed direction $\theta$ since the dual picture can often be much easier to analyse. For instance, a question raised in \cite{BSS22}, was to compute the Hausdorff dimension of $\NU_0^\theta$, the set of points with non-unique $\theta$-directed downward semi-infinite geodesics. Though it is unclear how to compute this dimension directly, one can in fact use duality to convert the above problem into computing the dimension of a geodesic, which is easy to show to be $4/3$ owing to the $2/3-$ H\"older continuity of geodesics (see Propositions \ref{prop:3}, \ref{prop:4}). We now state the above result along with an analogous result for the set of points admitting three distinct geodesics in a fixed direction, and these form the new contributions of this paper regarding the fractal geometry of the directed landscape. %
\begin{theorem}
  \label{thm:6}
  For any fixed direction $\theta\in \RR$ and the coupling $(\cL_\theta,\wcL_\theta)$ from Theorem \ref{thm:20}, the following statement holds almost surely. We have $\NU_0^\theta(\cL_\theta)=\cI^\theta_\uparrow(\cL_\theta)=\cT^\theta_\uparrow(\wcL_\theta)$, and this set has Hausdorff dimension $4/3$. Further, the set of points $p$ admitting three distinct geodesics $\Gamma^\theta_p$ is equal to $\mathrm{Tri}(\cI^\theta_\uparrow (\cL_\theta))=\mathrm{Tri}(\cT^\theta_\uparrow(\wcL_\theta))$ and is countable.
\end{theorem}

\paragraph{\textbf{Notational comments}}
We will usually use $\0$ to denote the point $(0,0)\in \RR^2$ and will often use $\leb$ to denote the Lebesgue measure on $\RR$. Throughout the paper, we will have upward and downward versions of different objects like interfaces, geodesics, geodesic trees and interface portraits. We will often only state results for one of the two versions, and it will be understood that corresponding results hold for the other version as well. As mentioned earlier, we will refer to a continuous function $\eta\colon [s,t]\rightarrow \RR$ as a path and additionally call it a semi-infinite or bi-infinite path if one or both of $s,t$ are $\pm \infty$. For any such path $\eta$, we will use $\inte (\eta)$ to denote the path $\eta\lvert_{(s,t)}$ as well as its graph, by which we mean the set $\{(\eta(r),r): r\in (s,t)\}\subseteq \RR^2$. %
Throughout the paper, we will often work with multiple paths emanating (starting) out of the same point (usually geodesics or interfaces), and to denote the left-most and right-most such paths, we respectively underline and overline the corresponding variable, for e.g., $\underline{\gamma}_{(x,s)}^{(y,t)}$ will denote the left-most geodesic from $(x,s)$ to $(y,t)$. Also, as commonly used, we will use the word `countable' for sets having the cardinality of $\NN$ and will use the phrase `at most countable' for sets with cardinality at most that of $\NN$.  
\paragraph{\textbf{Acknowledgements}}
The author thanks Riddhipratim Basu and Shirshendu Ganguly for the comments and acknowledges the support of the Institute for Advanced Study and the NSF grants DMS-1712862 and DMS-1953945. The author also acknowledges the International Centre for Theoretical Sciences, Bangalore for their hospitality since much of the work was completed during a visit there. The author thanks B{\'a}lint Vir{\'a}g and an anonymous referee for pointing out the connection of this paper to the Brownian web. Finally, the author thanks the anonymous referees for their very detailed and helpful comments and especially the suggestion to use the definition of interface portraits in \eqref{eq:15} instead of the more technical definition used earlier-- these comments have greatly improved the paper.
\section{The graph structure of semi-infinite geodesics and interfaces}
\label{sec:tree}
\subsection{Required results from the literature}
\label{sec:import}
To begin, we state some known results about the directed landscape %
which will be useful for the present work. Being a natural scaling limit of discrete models, the directed landscape has several symmetries, and we record these in the following proposition.

\begin{proposition}[{\cite[Lemma 10.2]{DOV18} and \cite[Proposition 1.23]{DV21}}]
  \label{prop:symm}
  When considered as random continuous functions from $\RR_\uparrow ^4\rightarrow \RR$, the following distributional equalities hold.
  \begin{itemize}
  \item KPZ $1\colon 2 \colon 3$ scaling: For any $q>0$, $\cL(x,s;y,t)\stackrel{d}{=}q\cL(q^{-2}x,q^{-3}s;q^{-2}y,q^{-3}t)$.
  \item Translational invariance: For any $x_0,s_0\in \RR$, $\cL(x,s;y,t)\stackrel{d}{=}\cL(x+x_0,s+s_0;y+x_0,t+s_0)$.
  \item Flip invariance: $\cL(x,s;y,t)\stackrel{d}{=}\cL(-x,s;-y,t)\stackrel{d}{=}\cL(y,-t;x,-s)$.
  \item Skew invariance: For any $\theta\in \RR$,
    \begin{displaymath}
      \cL(x,s;y,t)\stackrel{d}{=}\cL(x+\theta s,s;y+\theta t, t) +\theta^2(t-s) +2\theta(y-x).
    \end{displaymath}
  \end{itemize} 
  
\end{proposition}
In fact, the skew invariance discussed above transforms geodesic trees and interface portraits in a simple manner as we shall see now.

\begin{lemma}
  \label{lem:skew}
  Given the landscape $\cL$ and a fixed $\theta\in \RR$, define $\cL^{\mathrm{sk}}_\theta\colon \RR_\uparrow^4\rightarrow \RR$ by $\cL_\theta^{\mathrm{sk}}(x,s;y,t)=\cL(x+\theta s,s;y+\theta t, t) +\theta^2(t-s) +2\theta(y-x)$ which, by Proposition \ref{prop:symm}, is also a directed landscape. Then, with $\pi_\theta\colon \RR^2\rightarrow \RR^2$ denoting the skew transformation $\pi_\theta(x,s)= (x+\theta s,s)$, we almost surely have $\cT^\theta_\downarrow(\cL)=\pi_\theta(\cT_\downarrow(\cL^\mathrm{sk}_\theta))$ and $\cI^\theta_\uparrow(\cL)=\pi_\theta(\cI_\uparrow(\cL^\mathrm{sk}_\theta))$.
\end{lemma}
\begin{proof}
  It suffices to show the equality for the geodesic trees as then the corresponding equality for the interface portraits would immediately follow by definition.  Now, given a path $\eta\colon [s_1,t_1]\rightarrow \RR$, let $\pi_\theta\eta$ be the skew-transformed path defined by $\pi_\theta\eta(s)=\eta(s)+\theta s$. We now note that to show the a.s.\ equality $\cT^\theta_\downarrow(\cL)=\pi_\theta(\cT_\downarrow(\cL^\mathrm{sk}_\theta))$, it suffices to show that almost surely, for any $(x',s';y',t')\in \RR_\uparrow^4$ and any $\cL_\theta^\mathrm{sk}$-geodesic $\gamma$ from $(x',s')$ to $(y',t')$, the path $\pi_\theta \gamma$ is an $\cL$-geodesic, and we now proceed to show this.

  With $\ell_{\cL}$ and $\ell_{\cL_\theta^\mathrm{sk}}$ being used to denote the lengths in the respective landscapes as defined in \eqref{eq:37}, note that for any path $\eta$ from $(x',s')$ to $(y',t')$, we have
  \begin{equation}
    \label{eq:2}
    \ell_{\cL_\theta^\mathrm{sk}}(\eta)=\ell_{\cL}(\pi_\theta\eta)+\theta^2(t'-s')+2\theta(y'-x').
  \end{equation}
 Since the term $\theta^2(t'-s')+2\theta(y'-x')$ has no dependence on the path $\eta$, we immediately obtain that if $\gamma$ maximizes the $\ell_{\cL_\theta^\mathrm{sk}}$ length among paths from $(x',s')$ to $(y',t')$, then $\ell_{\cL}(\pi_\theta\gamma)$ is the maximum over the $\ell_{\cL}$ lengths of paths from $(x'+\theta s',s')$ to $(y'+\theta t',t)$. This completes the proof.
\end{proof}

Later, we shall need to compute the Hausdorff dimension of semi-infinite geodesics, and for the upper bound, we shall use the following result on the H\"older continuity of geodesics.
\begin{proposition}[{\cite[Proposition 12.3]{DOV18}}]
  \label{prop:3}
 For any fixed point $(x,s;y,t)\in \RR_\uparrow^4$ and any fixed $\varepsilon>0$, the geodesic $\gamma_{(x,s)}^{(y,t)}$ is a.s.\ H\"older $2/3-\varepsilon$. %
\end{proposition}
On the other hand, for lower bounding the dimension of semi-infinite geodesics, we will require the following occupation measure estimate.
\begin{proposition}[{\cite[Lemma 2.7]{GZ22}}]
  \label{prop:4}
  There exist constants $C,c$ such that the following is true. With $\leb$ denoting the Lebesgue measure on $\RR$, for any closed intervals $I\subseteq \RR$ and $J\subseteq (-\infty,0]$, and for all $M>0$, we have
  \begin{displaymath}
    \PP(\leb( \{t\in J: \Gamma_{\0}(t)\in I\})> M \leb(I)\leb(J)^{1/3})<Ce^{-cM}.
  \end{displaymath}
\end{proposition}
We now state a bound on the transversal fluctuation of semi-infinite geodesics followed by a lemma which is an easy consequence.
\begin{proposition}[{\cite[Theorem 3.2]{RV21}}]
  \label{prop:1.1}
  There exists a random constant $N$ satisfying $\PP(N\geq t)\leq Ce^{-ct^3}$ for some positive constants $C,c$ and all $t>0$ such that for all $s>0$, 
  \begin{displaymath}
    |\Gamma_\0(-s)|\leq Ns^{2/3} (1+\log^{1/3}(|\log s|)). 
  \end{displaymath}
\end{proposition}
\begin{lemma}
  \label{lem:3}
For any fixed point $p=(y_0,t_0)\in \RR^2$ and any fixed $s_0<t_0$, we almost surely have $\lim_{n\rightarrow \infty}(\min_{s\in [s_0,t_0]}\Gamma_{p+(n,0)}(s))=\infty$ and $\lim_{n\rightarrow \infty}(\max_{s\in [s_0,t_0]}\Gamma_{p-(n,0)}(s))=-\infty$.
\end{lemma}
\begin{proof}
  By the flip symmetry of the directed landscape (Proposition \ref{prop:symm}), it suffices to show the former equality. By applying Proposition \ref{prop:1.1}, there are constants $C,c$ such that, $\PP(\min_{s\in [s_0,t_0]}\Gamma_{p+(n,0)}(s)<y_0+n/2)\leq Ce^{-cn^3}$ and thus by applying the Borel-Cantelli lemma, we obtain that almost surely, $\min_{s\in [s_0,t_0]}\Gamma_{p+(n,0)}(s)\rightarrow \infty$ as $n\rightarrow \infty$. This completes the proof.
\end{proof}
We shall also require a simultaneous bound on the transversal fluctuations of all finite geodesics, and we now state such a result.
\begin{proposition}[Lemma 3.11 in \cite{GZ22}]
  \label{prop:trans-simul}
  There is a positive random variable $S$ and constants $C,c>0$ such that the following holds. For $M>0$, we have $\PP(S>M)< Ce^{-cM^{9/4}(\log M)^{-4}}$. Moreover, for any $u=(x,s;y,t)\in \RR_{\uparrow}^4$, and any geodesic $\gamma_{(x,s)}^{(y,t)}$ and $(s+t)/2\leq r<t$,
  \begin{displaymath}
    \left|
      \gamma_{(x,s)}^{(y,t)}(r)-\frac{x(t-r)+y(r-s)}{t-s}
    \right|< S(t-r)^{2/3}\log^3
    \left(
      1+\frac{\|u\|_2}{t-r}
    \right),
  \end{displaymath}
where $\|u\|_2$ denotes the usual $\ell^2$ norm.  A similar bound holds when $s<r<(s+t)/2$ by symmetry.
\end{proposition}
We will often want to take the limit of a sequence of geodesics, and the following result from \cite{DSV22} will ensure that subsequential limits always exist and that the limits are always geodesics. However, due to minor technical issues regarding functions with possibly different domains of definition, we first need to clarify the precise notion of convergence of paths used here (and throughout the paper), and we do so now. For any $s<t\in \RR\cup\{\pm \infty\}$, sequences $\{s_n\}_{n\in \NN},\{t_n\}_{n\in \NN}\subseteq \RR\cup\{\pm \infty\}$ satisfying $\lim_{n\rightarrow \infty}s_n=s$, $\lim_{n\rightarrow \infty}t_n=t$, and paths $\gamma_n\colon [s_n,t_n]\rightarrow \RR$, we say that $\gamma_n$ converges to $\gamma$ uniformly if we have $\lim_{n\rightarrow \infty}\sup_{s'\in [s,t]\cap[s_n,t_n]}|\gamma_n(s')-\gamma(s')|=0$. Similarly, we say that $\gamma_n$ converges to $\gamma$ locally uniformly if for every bounded open interval $I\subseteq [s,t]$, we have $\lim_{n\rightarrow \infty}\sup_{s'\in I\cap[s_n,t_n]}|\gamma_n(s')-\gamma(s')|=0$. We are now ready to state the precompactness result.

\begin{proposition}[{\cite[Lemma 3.1]{DSV22}}] 
  \label{prop:9}
The following holds with probability $1$. For all points $u=(x,s;y,t)\in \RR^4_\uparrow$ and any sequence of points $u_n=(x_n,s_n;y_n,t_n)\in \RR^4_\uparrow$ converging to $u$, every sequence of geodesics $\gamma_{(x_n,s_n)}^{(y_n,t_n)}$ is precompact in the uniform metric, and every subsequential limit is a geodesic from $(x,s)$ to $(y,t)$. 
\end{proposition}
Similar to the above, we have the following precompactness result for semi-infinite geodesics.
\begin{lemma}
  \label{prop:12}
 The following holds with probability $1$. For all $\theta\in \RR$, all $p=(y,t)\in \RR^2$ and any sequence of points $p_n\rightarrow p$, every sequence of geodesics $\Gamma^\theta_{p_n}$ is precompact in the locally uniform topology and every subsequential limit is a geodesic $\Gamma^\theta_{p}$.
\end{lemma}
\begin{proof}
  First, we show that the sequence $\Gamma^\theta_{p_n}(t-1)$ is a.s.\ bounded, and for this part, we just cite results from \cite{BSS22}. Indeed, \cite[Theorem 6.5 (i)]{BSS22} states that all the geodesics $\underline{\Gamma}^\theta_{p_n}$ and $\overline{\Gamma}^\theta_{p_n}$ are in fact ``Busemann geodesics'' and then \cite[Theorem 5.9 (i)]{BSS22} stated for ``Busemann geodesics'' implies that the set comprised of points of the form $\underline{\Gamma}^\theta_{p_n}(t-1),\overline{\Gamma}^\theta_{p_n}(t-1)$ is a.s.\ bounded. As a result, the sequence $\Gamma^\theta_{p_n}(t-1)$ is a.s.\ bounded as well.
  
   Now, we let $M$ be a random positive integer such that $|\Gamma^\theta_{p_n}(t-1)|<M$ for all large enough $n$ and consider the geodesics $\overline{\Gamma}_{(M,t-1)}^\theta,\underline{\Gamma}_{(-M,t-1)}^\theta$. Due to the left-most and right-most nature of the above two geodesics, it is not difficult to see that for all large enough $n$ and all $t'\leq t-1$, we have
  \begin{equation}
    \label{eq:45}
    \underline{\Gamma}_{(-M,t-1)}^\theta(t')\leq\Gamma^\theta_{p_n}(t')\leq \overline{\Gamma}_{(M,t-1)}^\theta(t')
  \end{equation}
  and thus, since we know that $\lim_{t'\rightarrow -\infty}\underline{\Gamma}_{(-M,t-1)}^\theta(t')/t'=\lim_{t'\rightarrow -\infty}\overline{\Gamma}_{(M,t-1)}^\theta(t')/t'=\theta$, we obtain that any subsequential limit of the geodesics $\Gamma_{p_n}^\theta$ must be downward $\theta$-directed as well.

Thus, in view of the above and Proposition \ref{prop:2}, it suffices to show that the geodesics $\Gamma_{p_n}^\theta$ are precompact in the locally uniform topology.  For this, note that by \eqref{eq:45}, for any $N\in \NN$, the points $\Gamma_{p_n}^\theta(t-N)$ are bounded in $n$ and thus they converge along a subsequence to a point $x_N$. By using Proposition \ref{prop:9} with $u=(x_N,t-N;y,t)$, we obtain that the geodesics $\Gamma_{p_n}^\theta\lvert_{[t-N,t_n]}$ are precompact in the uniform topology. Since $N$ is arbitrary, we obtain that the geodesics $\Gamma_{p_n}^\theta$ are precompact in the locally uniform topology, and this completes the proof.
\end{proof}

The following result states that if a sequence of geodesics converges to another geodesic uniformly, then it must also converge in the ``overlap sense''.
\begin{proposition}[{\cite[Lemma B.12]{BSS22},\cite[Lemma 3.3]{DSV22}}]
  \label{prop:1.2}
  The following holds with probability $1$. Let $u=(x,s;y,t)\in \RR^4_\uparrow$ and $u_n=(x_n,s_n;y_n,t_n)\in \RR^4_\uparrow$ be a sequence converging to $u$. Suppose that either there is a unique geodesic $\gamma_{(x,s)}^{(y,t)}$ or the geodesics $\gamma_{(x_n,s_n)}^{(y_n,t_n)}$ are unique and converge uniformly to a geodesic $\gamma_{(x,s)}^{(y,t)}$. Then the geodesics $\gamma_{(x_n,s_n)}^{(y_n,t_n)}$ converge to $\gamma_{(x,s)}^{(y,t)}$ in the overlap sense, by which we mean that the set
  \begin{equation}
    \label{eq:43}
    \left\{s'\in [s_n,t_n]\cap [s,t]\colon \gamma_{(x_n,s_n)}^{(y_n,t_n)}(s')=\gamma_{(x,s)}^{(y,t)}(s')\right\}
  \end{equation}
  is an interval whose endpoints converge to $s$ and $t$. 
\end{proposition}
Though we shall not require its full strength, the following result from \cite{Bha22} concerning atypical points on geodesics will be useful to us.
\begin{proposition}[{\cite[Theorem 4]{Bha22}}]
  \label{prop:8}
Call a point $p\in \RR^2$ an atypical star if it admits two disjoint geodesics to points $q_1,q_2$ both of which are on the same side of $p$ with respect to the time co-ordinate.  For any fixed $u=(x_1,s_1;y_1,t_1)\in \RR_\uparrow^4$, consider the set of $t\in (s_1,t_1)$ such that the point $(\gamma_{(x_1,s_1)}^{(y_1,t_1)}(t),t)$ is an atypical star. Then the above set almost surely has Hausdorff dimension $1/3$.
\end{proposition}
The following result has already been discussed earlier, but since it is important to this work, we state it as a proposition.
\begin{proposition}[{\cite[Theorem 2.5]{BSS22}}]
  \label{prop:11}
  There is a random countable set $\Xi_\downarrow\subseteq \RR$ such that almost surely, for any $\theta\in \Xi_\downarrow^c$, any two downward $\theta$-directed semi-infinite geodesics $\Gamma^\theta_p,\Gamma^\theta_q$ for points $p,q\in \RR^2$ eventually coalesce. A similar result holds for upward semi-infinite geodesics and the corresponding non-uniqueness set $\Xi_\uparrow$.
\end{proposition}

Finally, we state a result from \cite{BSS22} (see also \cite[Theorem 3.9]{RV21}) regarding the cardinality of the time-slices of the sets $\NU_0^\theta$.
\begin{proposition}[{\cite[Theorem 6.1]{BSS22}}]
  \label{lem:6}
  The following holds with probability $1$. For all $\theta \in \Xi_\downarrow^c$ and all $t>0$, the set $\NU_0^\theta \cap \{(x,s)\in \RR^2:s=t\}$ is countable.
\end{proposition}
\subsection{Theorem \ref{thm:2}: Proof and consequences}
\label{sec:thm2pf}
The first aim of this section is to prove Theorem \ref{thm:2}, and this will be the crucial ingredient in showing that semi-infinite geodesics in any uniqueness direction form a tree. We begin with a few preparatory lemmas.

\begin{figure}
  \centering
  \includegraphics[width=0.2\textwidth]{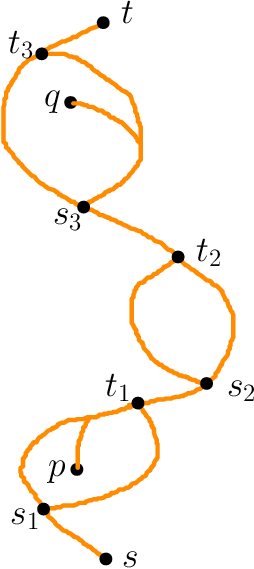}
  \caption{Proof of Lemma \ref{lem:59}: If $\eta$ has $(s_1;t_1),(s_2;t_2),(s_3;t_3)$ bubbles, then due to topological considerations, there must be exist at least two geodesics between $p$ and $q$, one using the left side of the $(s_2;t_2)$ bubble and the other using the right side. However, since the points $p,q$ are chosen to be rational, this is not possible.}
  \label{fig:three-bub}
\end{figure}
For a geodesic $\eta\colon[s,t]\rightarrow \RR$ and points $s',t'$ with $s<s'<t'<t$, we say that $\eta$ has an $(s';t')$ bubble if there exists a geodesic $\gamma_{(\eta(s'),s')}^{(\eta(t'),t')}$ which is disjoint from $\eta$ except at its endpoints.
\begin{lemma}
  \label{lem:59}
  Almost surely, there does not exist any $s<t$ and a geodesic $\eta\colon[s,t]\rightarrow \RR$ such that there exist points $s<s_1<t_1<s_2<t_2<s_3<t_3<t$ with $\eta$ having $(s_1;t_1),(s_2;t_2),(s_3;t_3)$ bubbles. 
\end{lemma}
\begin{proof}
    If such a geodesic were to exist, then we could take a rational point (see Figure \ref{fig:three-bub}) $p$ lying inside the $(s_1;t_1)$ bubble and a rational point $q$ lying inside the $(s_3;t_3)$ bubble and we would have multiple geodesics from $p$ to $q$, but this is impossible since rational points have unique geodesics between them (see \cite[Theorem 12.1]{DOV18}).
\end{proof}
\begin{lemma}
  \label{lem:58}
Recall the notion of an atypical star from Proposition \ref{prop:8}.  Almost surely, for any geodesic $\eta\colon[s',t']\rightarrow\RR$ and any $(s,t)\subseteq (s',t')$, there exists $r\in (s,t)$ such that $(\eta(r),r)$ is not an atypical star.
\end{lemma}
\begin{proof}
Choose points $s_j,t_j\in(s,t)$ for $j\in \{1,2,3\}$ as stated in Lemma \ref{lem:59} and let $j_0$ be the index for which $\eta$ does not have a $(s_{j_0};t_{j_0})$ bubble. Let $p_n,q_n$ be rational points approximating $p=(\eta(s_{j_0}),s_{j_0})$ and $q=(\eta(t_{j_0}),t_{j_0})$ respectively and let $\gamma_n$ denote the unique geodesic from $p_n$ to $q_n$. By Proposition \ref{prop:9}, there is a geodesic $\gamma$ between $p,q$ such that along a subsequence $\{n_i\}$, the geodesics $\gamma_n$ converge to $\gamma$ uniformly. Since $\eta$ does not have a $(s_{j_0};t_{j_0})$ bubble, we know that there exists a $\delta>0$ such that either $\eta\lvert_{[s_{j_0},s_{j_0}+\delta]}=\gamma\lvert_{[s_{j_0},s_{j_0}+\delta]}$ or $\eta\lvert_{[t_{j_0}-\delta,t_{j_0}]}=\gamma\lvert_{[t_{j_0}-\delta,t_{j_0}]}$. By Proposition \ref{prop:1.2}, we also know that $\gamma_{n_i}$ converges to $\gamma$ in the overlap sense and this along with the above implies that there must exist an interval $I \subseteq (s_{j_0},t_{j_0})\subseteq (s,t)$ such that $\eta\lvert_I=\gamma\lvert_I=\gamma_n\lvert_I$ for all large $n$, and we now fix such an $n$ and such an interval $I$. By Proposition \ref{prop:8}, we know that the set of $r\in I$ such that $(\gamma_n(r),r)$ is an atypical star a.s.\ has Hausdorff dimension $1/3$, and thus there must exist at least one $r\in I$ such that $(\eta(r),r)$ is not an atypical star.
\end{proof}

\begin{figure}
  \centering
  \includegraphics[width=0.3\textwidth]{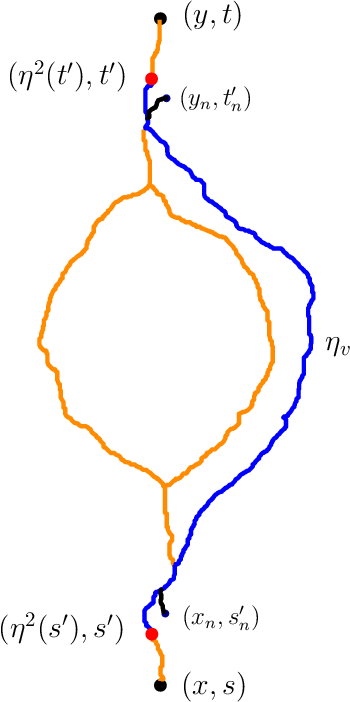}
  \caption{Proof of Theorem \ref{thm:2}: By using that points which are not atypical stars exist on any segment of any geodesic, first find $s',t'$ such that $(\eta^2(s'),s'),(\eta^2(t'),t')$ are not atypical stars, thereby implying that all geodesics $\eta_v$ corresponding to $(\eta^2(s'),s';\eta^2(t'),t')$ follow the orange geodesics initially. Now if the rational approximants $(x_n,s_n'),(y_n,t_n')$ are close enough to $(\eta^2(s'),s'),(\eta^2(t'),t')$, then due to coalescence, the orange bubble will lead to multiple geodesics for the rational point $(x_n,s_n';y_n,t_n')$ which cannot happen.}
  \label{fig:bubble}
\end{figure}

With the above preparatory lemmas at hand, we are ready to complete the proof of Theorem \ref{thm:2}. 
\begin{proof}[Proof of Theorem \ref{thm:2}]
  Suppose (with the aim of eventual contradiction) that there exists a $u=(x,s;y,t)\in \RR_\uparrow^4$ and two distinct geodesics $\eta^1,\eta^2$ from $(x,s)$ to $(y,t)$ with the property that for a small enough $\delta>0$, $\eta^1(r)=\eta^2(r)$ for all $r\in(s,s+\delta)\cup (t-\delta,t)$. We claim that there exist $s'\in (s,s+\delta),t'\in (t-\delta,t)$ and a $\delta'>0$ such that the point $v=(\eta^2(s'),s';\eta^2(t'),t')$ has the property that for any geodesic $\eta_v$ associated to $v$, we have $\eta_v(r)=\eta^1(r)=\eta^2(r)$ for all $r\in (s',s'+\delta')\cup (t'-\delta',t')$. Indeed, this is true because by Lemma \ref{lem:58}, we can find such $s'\in (s,s+\delta),t'\in (t-\delta,t)$ such that none of $(\eta^2(s'),s')$ and $(\eta^2(t'),t')$ are atypical stars. For the remainder of the argument, we fix such a choice of $s',t',\delta'$ and the corresponding $v$.%

 With the aim of eventually exhibiting a rational point in $\RR^4_\uparrow$ with multiple geodesics corresponding to it, we choose $u_n=(x_n,s_n';y_n,t_n')\in \RR_\uparrow^4$ to be a rational sequence converging to $v$ and let $\psi_n$ denote the unique geodesics corresponding to $u_n$. Note that by Proposition \ref{prop:9}, there exists a subsequence  $\{n_i\}$ along which $\psi_n$ converges to a geodesic $\eta_v$ corresponding to $v$. In fact, if we define the set $I_n$ by $I_n=\{r\in \RR:\psi_n(r)=\eta_v(r)\}$, then by Proposition \ref{prop:1.2}, $I_n$ is an interval which increases to $[s',t']$ as $n\rightarrow \infty$. Now, we choose an $n$ large enough so that $(s',s'+\delta')\cap I_n$ and $(t'-\delta',t')\cap I_n$ are both non-empty, and then we observe that there are at least two geodesics (see Figure \ref{fig:bubble}) between $(x_n,s_n')$ and $(y_n,t_n')$, and this yields a contradiction as all rationals simultaneously have unique geodesics between them.

\end{proof}
We now use Theorem \ref{thm:2} to prove Corollary \ref{thm:3}.
\begin{proof}[Proof of Corollary \ref{thm:3}]
  Fix $\theta\in \Xi_\downarrow^c$. By definition, it is clear that $\NU^\theta_1\subseteq \NU^\theta_0$. For the reverse implication, first note that by Proposition \ref{prop:11}, for any $(x,s)\in \NU^\theta_0$, we have $\underline{\Gamma}_{(x,s)}^\theta(t)= \overline{\Gamma}_{(x,s)}^\theta(t)$ for all negative enough $t$, and we fix such a choice of $t$. Now, if $(x,s)\notin \NU^\theta_1$, then we consider $u=(\Gamma^\theta_{(x,s)}(t-1),t-1;x,s)\in \RR_\uparrow^4$, and we note that such a point cannot exist by Theorem \ref{thm:2}. Thus $(x,s)\in \NU^\theta_1$, thereby establishing that $\NU^\theta_0\subseteq \NU^\theta_1$ and subsequently $\NU^\theta_0=\NU^\theta_1$.
\end{proof}

We spend the rest of this section obtaining some useful consequences of Theorem \ref{thm:2}. The first consequence concerns the overlap approximation of semi-infinite geodesics by those started from nearby points.

\begin{lemma}
  \label{lem:47}
  The following hold with probability $1$ simultaneously for all $\theta \in \Xi_\downarrow^c$.
  \begin{enumerate}
  \item For any $\varepsilon>0$ and any point $p=(y,t)$ with a unique geodesic $\Gamma_p^\theta$, there exists a neighbourhood $U_\varepsilon$ such that all points $q\in U_\varepsilon$ and all geodesics $\Gamma^\theta_q$ satisfy $\Gamma^\theta_q\lvert_{(-\infty,t-\varepsilon)}=\Gamma^\theta_p\lvert_{(-\infty,t-\varepsilon)}$.
  \item For any point $p=(y,t)$ with a not necessarily unique geodesic $\Gamma^\theta_p$ and any $\varepsilon>0$, there exists a non-empty open set $V_\varepsilon$ such that for all $q\in V_\varepsilon$ and all geodesics $\Gamma_q^\theta$, we have $\Gamma^\theta_p\lvert_{(-\infty,t-\varepsilon)}\subseteq \inte(\Gamma^\theta_q)$.
  \end{enumerate}
\end{lemma}

\begin{proof}
  By applying Lemma \ref{prop:12}, we can choose a neighbourhood $U_\varepsilon$ of $p$, so as to make the quantity $\sup_{q\in U_\varepsilon,s\in (t-\varepsilon,t-\varepsilon/2)}|\Gamma^\theta_p(s)-\Gamma^\theta_q(s)|$ arbitrarily small. In conjunction with Proposition \ref{prop:1.2}, we can thus ensure that $\Gamma^\theta_p,\Gamma^\theta_q$ meet in the time interval %
  $(t-\varepsilon,t)$ for all $u\in U_\varepsilon$ and geodesics $\Gamma_q^\theta$. Since $\Gamma_p^\theta$ is unique, the two geodesics must coalesce once they meet, and this completes the proof of the first statement. %
  
  To obtain the second statement, we note that by Theorem \ref{thm:2}, for any point $p$ with a not necessarily unique geodesic $\Gamma_p^\theta$, the geodesic $\Gamma_p^\theta\lvert_{(-\infty,t-\varepsilon/2]}$ is the unique $\theta$-directed downward semi-infinite geodesic emanating from its end point. We can now apply the first statement of the lemma, and this finishes the proof of the second statement.
\end{proof}

As an immediate consequence of the second statement above, we have the following result.
\begin{lemma}
  \label{lem:8}
  For any fixed $\theta\in \RR$, we have the a.s.\ equality $\cT^\theta_\downarrow = \bigcup_{q\in \QQ^2}\inte(\Gamma^\theta_{q})$.
\end{lemma}

Theorem \ref{thm:2} can be used to obtain another interesting consequence. The result \cite[Corollary 3.5]{DSV22} states that that the union of interiors of all left-most and right-most geodesics is equal to the union of interiors of rational geodesics. By using Theorem \ref{thm:2}, we can upgrade the above to a corresponding statement for the union of the interiors of all possible geodesics, instead of just the left-most and right-most ones. %
\begin{proposition}
  \label{lem:50}
  Borrowing terminology from the Brownian geometry literature \cite{MQ20, AKM17}, define the geodesic frame $\cW$ to be the union of the interiors of all geodesics in the directed landscape. Then 
  \begin{equation}
    \label{eq:46}
  \cW=\bigcup_{(x,s;y,t)\in \QQ_\uparrow^4}\inte(\gamma_{(x,s)}^{(y,t)}).  
  \end{equation}
   Further, for any fixed point $p\in \RR^2$, $p\notin \cW$ with probability $1$, and as a consequence, $\cW$ is a.s.\ a strict subset of $\RR^2$.
\end{proposition}
\begin{proof}
We begin by proving \eqref{eq:46}. That the latter set is a subset of the former is clear, and it suffices to show the reverse inclusion. For any $(x',s';y',t')\in \RR_\uparrow^4$, consider a sequence of intervals $[s'_n,t'_n]$ increasing to $[s',t']$. As a consequence of Theorem \ref{thm:2}, $\gamma_{(x',s')}^{(y',t')}\lvert_{[s'_n,t'_n]}$ is the unique geodesic between its endpoints, and thus by Proposition \ref{prop:1.2}, %
  we obtain that $\inte(\gamma_{(x',s')}^{(y',t')}\lvert_{[s_n',t_n']})\subseteq \bigcup_{(x,s;y,t)\in \QQ_\uparrow^4}\inte(\gamma_{(x,s)}^{(y,t)})$. Using that $\inte(\gamma_{(x',s')}^{(y',t')})=\bigcup_{n\in \NN}\inte(\gamma_{(x',s')}^{(y',t')}\lvert_{[s_n',t_n']})$, we obtain that $\inte(\gamma_{(x',s')}^{(y',t')})\subseteq \bigcup_{(x,s;y,t)\in \QQ_\uparrow^4}\inte(\gamma_{(x,s)}^{(y,t)})$, and this shows that $\cW=\bigcup_{(x,s;y,t)\in \QQ_\uparrow^4}\inte(\gamma_{(x,s)}^{(y,t)})$. %

We now come to the second statement. By the first statement along with the translational symmetry of the directed landscape (Proposition \ref{prop:symm}), it suffices to show that for any fixed $(x,s;y,t)\in \RR_\uparrow^4$ with $s<0<t$, $\0$ does not lie on the geodesic $\gamma_{(x,s)}^{(y,t)}$. However, this is immediate by observing that $\gamma_{(x,s)}^{(y,t)}(0)$ is the location of unique maximizer of the profile $x'\mapsto (\cL(x,s;x',0)+\cL(x',0;y,t))$, where the two summands are independent parabolic $\mathrm{Airy}_2$ processes (see \cite[Definition 10.1 (I)]{DOV18}). Indeed, by the locally Brownian nature of the parabolic $\mathrm{Airy}_2$ process \cite{CH14} (more specifically, see the proof of \cite[Theorem 4.3]{CH14}), it is easy to see that there is zero probability of the above maximum being achieved at $x'=0$. This completes the proof.
\end{proof}
Now, as a simple consequence of Proposition \ref{lem:50}, we obtain the following non-triviality of the objects $\cT_\downarrow^\theta$ for all $\theta\in \Xi_\downarrow^c$.
\begin{lemma}
  \label{lem:1}
  For any fixed point $p\in \RR^2$, we almost surely have $p\notin \cT_\downarrow^\theta$ simultaneously for all $\theta \in \Xi_\downarrow^c$. In particular, almost surely, $\cT^\theta_\downarrow$ is a strict subset of $\RR^2$ simultaneously for all $\theta \in \Xi_\downarrow^c$.
\end{lemma}
\begin{proof}
  First note that almost surely, for all $\theta \in \Xi_\downarrow^c$, we have $\cT_\downarrow^\theta\subseteq \cW$ and then apply Proposition \ref{lem:50}.
\end{proof}
Now, we state a result demonstrating that when traversing between bounded and disjoint regions of space, geodesics need to pass via finitely many ``highways''; this was proved for exponential LPP in \cite{BHS22}, and this result will be important to obtain the duality later in the paper. %
\begin{lemma}
  \label{lem:45}
  For fixed values of $a<b,s<s'<t'<t$, let $\cC$ denote the collection of all paths obtained as $\gamma_{(x,s)}^{(y,t)}\lvert_{[s',t']}$ for some $x,y\in [a,b]$. Then the set $\cC$ is a.s.\ finite. Further, there exist constants $C,c>0$ depending on $a,b,s,t,s',t'$ such that for all $M>0$, we have
  \begin{equation}
    \label{eq:40}
    \PP(|\cC|\geq M)\leq Ce^{-cM^{1/384}}.
  \end{equation}
\end{lemma}
 In order to facilitate the proof of the above, we now state a tail bound, obtained in \cite{GZ22} by using the results from \cite{BHS22}, on the number of points on a horizontal line which are hit by geodesics between macroscopic regions which are additionally unique between their endpoints.
\begin{proposition}[{\cite[Lemma 3.12]{GZ22})}]
  \label{prop:finhits}
  Fix real numbers $a<b,s<r<t$. Then there exist constants $C',c'$ depending only on $(r-s)/(t-s)$ such that for all $M>(b-a)^3/(t-s)^2$,
  \begin{displaymath}
    \PP(|\{\gamma_{(x,s)}^{(y,t)}(r):x,y\in [a,b], \gamma_{(x,s)}^{(y,t)} \textrm{ is unique}\}|>M)\leq C'e^{-c'M^{1/384}}.
  \end{displaymath}
\end{proposition}
We note that \cite{GZ22} states the above result with $s=0$ and $t=1$ and the version of the result stated above can be obtained by combining this with the KPZ and translational symmetries of the directed landscape (Proposition \ref{prop:symm}). We now combine Proposition \ref{prop:finhits} with Theorem \ref{thm:2} to complete the proof of Lemma \ref{lem:45}.
\begin{proof}[Proof of Lemma \ref{lem:45}]
We begin by setting up a good event $\cE_M$ on which transversal fluctuations of geodesics are controlled and then subsequently return to the proof. By Proposition \ref{prop:trans-simul}, we know that for some constants $C_1,c_1$, and on an event $\cE_M$ satisfying
  \begin{equation}
    \label{eq:51}
    \PP(\cE_M)\geq 1-C_1e^{-c_1M^{9/16}(\log M)^{-4}},
  \end{equation}
  we have
  \begin{equation}
    \label{eq:52}
    \sup_{x,y\in[a,b],r\in\{(s+s')/2,(t'+t)/2\}}|\gamma_{(x,s)}^{(y,t)}(r)|\leq M^{1/4}.
  \end{equation}
  With the event $\cE_M$ at hand, we now return to the proof. As an immediate consequence of Theorem \ref{thm:2}, observe that for geodesic $\gamma_{(x,s)}^{(y,t)}$ with $x,y\in [a,b]$, the values $\gamma_{(x,s)}^{(y,t)}(s')$ and $\gamma_{(x,s)}^{(y,t)}(t')$ entirely determine the path $\gamma_{(x,s)}^{(y,t)}\lvert_{[s',t']}$, and thus it suffices to bound the number of attainable pairs $(\gamma_{(x,s)}^{(y,t)}(s'),\gamma_{(x,s)}^{(y,t)}(t'))$. Secondly, again as a consequence of Theorem \ref{thm:2}, we note that for any geodesic $\gamma_{(x,s)}^{(y,t)}$ with $x,y\in [a,b]$, the path $\gamma_{(x,s)}^{(y,t)}\lvert_{[(s+s')/2,(t'+t)/2]}$ is the unique geodesic between its endpoints.
As a result of the above observations and the bound \eqref{eq:51}, it suffices to show that for some constants $C_2,c_2$ depending on $a,b,s,t,s',t'$, we have the bounds
  \begin{align}
    \label{eq:53}
    \PP(|\{\gamma_{(x,(s+s')/2)}^{(y,(t'+t)/2)}(s')\colon |x|,|y|\leq M^{1/4}, \gamma_{(x,(s+s')/2)}^{(y,(t'+t)/2)} \textrm{ is unique}\}|>M/2) &\leq C_2e^{-c_2M^{1/384}},\nonumber\\
    \PP(|\{\gamma_{(x,(s+s')/2)}^{(y,(t'+t)/2)}(t')\colon |x|,|y|\leq M^{1/4}, \gamma_{(x,(s+s')/2)}^{(y,(t'+t)/2)} \textrm{ is unique}\}|>M/2) &\leq C_2e^{-c_2M^{1/384}}.
  \end{align}
  Since $M/2> (2M^{1/4})^3/((t'+t-s-s')/2)^2$ for all large enough $M$, the above estimates follow immediately from Proposition \ref{prop:finhits}, and this completes the proof.
\end{proof}
We now use the finiteness in Lemma \ref{lem:45} to show that there are only countably many points at which geodesics coalesce.
\begin{lemma}
  \label{lem:44}
  Call $p\in \RR^2$ a confluence point if there exist two geodesics starting from some point $q\neq p$ which first split at $p$. Then almost surely, there are only countably many confluence points.
\end{lemma}

\begin{proof}
  For a rectangle $R_{a,b,s,t}=[a,b]\times [s,t]$ with $a,b,s,t\in \QQ$, let $\mathrm{Conf}(R_{a,b,s,t})$ denote the set consisting of confluence points obtained by all geodesics $\gamma_p^q,\gamma_p^r$ for all $p\in [a,b]\times \{s\}$ and $q,r\in [a,b]\times \{t\}$ along with the confluence points obtained by all geodesics $\gamma_q^p,\gamma_r^p$ for $q,r\in [a,b]\times \{s\}$ and $p\in [a,b]\times \{t\}$. Then for $n\in \NN$, by applying Lemma \ref{lem:45} with $a,b,s,t$ as above and with $s'=s+2^{-n},t'=t-2^{-n}$ and taking a union bound over all $n\in \NN$, we obtain that for any fixed $a,b,s,t$ as above, the set $\mathrm{Conf}(R_{a,b,s,t})$ is a.s.\ at most countable.

 Now, it is not difficult to see that the set of all confluence points is equal to $\bigcup\mathrm{Conf}(R_{a,b,s,t})$ where the union is taken over all rational rectangles. Indeed, to see this, suppose that for some $(q;p)=(x_q,s_q;x_p,s_p)\in \RR^4_\uparrow$, we have two geodesics $\eta_1,\eta_2$ which both start from $q$ and first split at $p$. In this case, for all $\delta$ small enough and $M$ large enough, all $a,b$ with $a<-M$ and $ b>M$, all $s\in (s_q,s_p)$ and all $t\in (s_p,s_p+\delta)$, we have $p\in \mathrm{Conf}(R_{a,b,s,t})$. This completes the proof of there almost surely being at most countably many confluence points.

We now show that there a.s.\ do exist infinitely many confluence points. To see this, we look at $\Gamma_{(n,0)}$ for all $n\in \NN$ and note that any two such geodesics create a confluence point at the first point where they meet each other, and it thus suffices to show that these special confluence points a.s.\ form an infinite set. By an application of Lemma \ref{lem:3}, for every $N\in \NN$, we almost surely have $\min_{s\in [-N,0]}\Gamma_{(n,0)}(s)\rightarrow \infty$ as $n\rightarrow \infty$, and it is not difficult to see that this implies that the above-mentioned set is a.s.\ infinite.
\end{proof}

We now come to the final result of this section, and for this, we need the language of directed forests which we now introduce. A downward directed forest $\fF$ is defined to be a collection of semi-infinite paths $\{\eta_i\}_{i\in I}$ for some indexing set $I$ with the additional property that any two such paths $\eta_i,\eta_j$ for $i\neq j$ coalesce when they first meet, where starting from the same point is not considered as ``meeting''. To be precise, this means that for any $i\neq j$ and paths $\eta_i\colon(-\infty,t_i]\rightarrow \RR$, if we have $\eta_i(t)=\eta_j(t)$ for some $t<\min(t_i,t_j)$, then we must have $\eta_i(s)=\eta_j(s)$ for all $s\leq t$.
 We denote any $\eta_i$ with $(\eta_i(t_i),t_i)=p\in \RR^2$ as $\fF_p$, with the understanding that there might be many possible choices for $\fF_p$ for a given $p\in \RR^2$. Note that we can similarly define upward directed forests as well. The following ordering property of directed forests is an immediate consequence of the definition.
\begin{lemma}
  \label{lem:7}
  For any downward directed forest $\fF$, no two paths $\fF_p$ and $\fF_q$ for any points $p,q\in \RR^2$ can ever cross, in the sense that we cannot find times $s< t$ such that one of $\fF_p(s)-\fF_q(s)$ and $\fF_p(t)-\fF_q(t)$ is strictly positive, while the other is strictly negative. 
\end{lemma}
\begin{proof}
  If some paths $\fF_p$, $\fF_q$ with the above property were to exist, then by continuity, there would have to exist a time $s'\in (s,t)$ for which $\fF_p(s')=\fF_p(s')$. However, by the definition of a downward directed forest, this would imply that we must also have $\fF_p(s)=\fF_q(s)$ which is a contradiction.
\end{proof}
 As an immediate consequence of Corollary \ref{thm:3}, we have the following result which, in view of Lemma \ref{lem:7}, we will often refer to as the ordering of semi-infinite geodesics.
\begin{lemma}
  \label{prop:2}
  Almost surely, simultaneously for all $\theta\in \Xi_\downarrow^c$, $\cT_\downarrow^\theta$ forms a downward directed forest.
\end{lemma}
We mention that when using directed forests, we will slightly abuse notation in the sense that apart from viewing $\fF$ as the collection $\{\fF_p\}_p$, we will also think of $\fF$ as the set $\bigcup_p \inte (\fF_p)\subseteq \RR^2$. For a concrete example, the geodesic tree $\cT_\downarrow$ is thought of both as the collection of all possible downward $0$-directed semi-infinite geodesics and as the set consisting of the union of interiors of all such geodesics.

\subsection{Basic properties of the interface portrait}
\label{sec:intface}
In this section, we record some basic properties of interfaces. %
To begin, we show that interface portraits $\cI_\downarrow^\theta$ form a directed forest simultaneously along all directions $\theta\in \Xi_\uparrow^c$. %

\begin{lemma}
  \label{lem:56}
 Almost surely, simultaneously for all $\theta\in \Xi_\uparrow^c$, $\cI^\theta_\downarrow$ forms a downward directed forest. 
\end{lemma}
\begin{proof}
 With the aim of eventually obtaining a contradiction, suppose that for some points $p=(y_1,t_1),q=(y_2,t_2)$ and some $\theta\in \Xi_\uparrow^c$, there exist interfaces $\Upsilon_p^\theta,\Upsilon_q^\theta$ and some $s<t<\min(t_1,t_2)$ such that $\Upsilon_p^\theta(t)=\Upsilon_q^\theta(t)$ but $\Upsilon_p^\theta(s)\neq \Upsilon_q^\theta(s)$. We now consider a point $q'=(x,s)\in \RR^2$ with $x$ lying strictly in between $\Upsilon_p^\theta(s)\neq \Upsilon_q^\theta(s)$. Now, by planarity, any geodesic $\Gamma_{q',\uparrow}^\theta$ must necessarily intersect one of $\inte(\Upsilon_p^\theta)$ or $\inte(\Upsilon_q^\theta)$, but this is a contradiction since, by definition, $\cI_\downarrow^\theta\cap \cT_\uparrow^\theta=\emptyset$.
\end{proof}

In order to discuss further properties, we need to introduce Busemann functions, a notion with roots in geometry \cite{Bus12} which will be useful for the arguments in this section as well as for the constructions in later sections. Originally introduced to first passage percolation in \cite{New95, Hof05}, Busemann functions enable us to interpret passage times to points ``infinitely far away'' in a given direction. For $p,q\in \RR^2$ and any $\theta\in \Xi_\downarrow^c$, we define
\begin{equation}
  \label{eq:35}
  \cB^\theta_\downarrow(p,q)=\cL(z;p)-\cL(z;q),
\end{equation}
where $z$ is the first meeting point of $\Gamma^\theta_p$ and $\Gamma^\theta_q$, and it is not difficult to see that the above is well-defined regardless of the choice of $\Gamma^\theta_p,\Gamma^\theta_q$, and is furthermore continuous \cite[Theorem 5.1]{BSS22} in both $p,q$. We note that the upward Busemann functions $\cB^\theta_\uparrow$ can be similarly defined as well. Now, by following \cite[Section 4.1]{RV21}, for any point $p=(y,t)\in \RR^2$, we introduce the competition function $d^\theta_{p,\uparrow}$ defined for $x\in \RR$ and $s<t$ by
\begin{equation}
  \label{eq:16}
   d_{p,\uparrow}^\theta(x,s)=\sup_{y'\geq y}
    \left\{
      \cB^\theta_{\uparrow}( (y',t),p)+\cL(x,s;y',t)
    \right\}-\sup_{y'\leq y}
    \left\{
      \cB^\theta_{\uparrow}( (y',t),p)+\cL(x,s;y',t)
    \right\}.
  \end{equation}
  As we shall see now in the following lemma consisting of an application of various results from \cite{RV21}, the utility of the above definition is that it yields natural choices of interfaces starting from all points $p$.
\begin{lemma}
  \label{lem:2}
  For each point $p=(y,t)\in \RR^2$, consider the functions $\underline{\Upsilon}_p^\theta, \overline{\Upsilon}_p^\theta\colon(-\infty,t]\rightarrow \RR$ defined by
  \begin{align*}
    \label{eq:18}
    \underline{\Upsilon}^\theta_p(s)&=\inf\{x\in \RR:d^\theta_{p,\uparrow}(x,s)\geq 0\},\nonumber\\
  \overline{\Upsilon}^\theta_p(s)&=\sup\{x\in \RR:d^\theta_{p,\uparrow}(x,s)\leq 0\}.
  \end{align*}
  Then almost surely, simultaneously for all $p\in \RR^2$ and $\theta\in \Xi_\uparrow^c$, $\underline{\Upsilon}^\theta_p$ and $\overline{\Upsilon}^\theta_p$ are downward $\theta$-directed interfaces starting at $p$. 
\end{lemma}
\begin{proof}
The proof consists of an application of a few results from \cite{RV21}, and this is the only location in the paper where these results are invoked. First, to bring ourselves to the setting of \cite{RV21}, we note that Proposition \ref{lem:6} implies that in the language of \cite[Section 3.4]{RV21}, the unique geodesic condition is satisfied simultaneously for all points $p=(y,t)\in \RR^2$, all $\theta\in \Xi_\uparrow^c$, and the initial conditions $y'\mapsto \cB_\uparrow^\theta( (y',t),p)$, where we remark that the geodesics $\Gamma_{q,\uparrow}^\theta$ for $q=(x,s)\in \RR^2$ with $s<t$ are viewed as geodesics to the initial condition given by $y'\mapsto \cB_\uparrow^\theta( (y',t),p)$ in the sense that
  \begin{equation}
    \label{eq:3}
    \Gamma_{q,\uparrow}^\theta(t)=\argmax_{y'\in \RR}\{\cL(q;y',t)+\cB_\uparrow^\theta((y',t),p)\}.
  \end{equation}
  Moreover, apart from the unique geodesic condition, Lemma \ref{prop:2} implies that the geodesic ordering condition from \cite[Lemma 3.8]{RV21} is satisfied simultaneously for all points $p=(y,t)\in \RR^2$, all $\theta\in \Xi_\uparrow^c$, and the initial conditions $y'\mapsto \cB_\uparrow^\theta( (y',t),p)$. As a result of these observations, the `good samples' assumption from \cite[Section 3.4]{RV21} holds a.s.\ for all these initial conditions simultaneously, and thus we are justified in using the results from \cite{RV21} simultaneously for all these initial conditions. We now move on to the proof.
   
To begin, we note that for any $\theta \in \Xi_\uparrow^c$, $d_{p,\uparrow}^\theta$ is a continuous as a function of $(x,s)$ due to the continuity of Busemann functions, and further, for any $s<t$, the map $x\mapsto d_{p,\uparrow}^\theta(x,s)$ is increasing \cite[Proposition 4.1]{RV21}. As a consequence, we obtain that $\underline{\Upsilon}_p^\theta(s)\leq \overline{\Upsilon}_p^\theta(s)$ for all $s\leq t$ and thus the former is indeed to the left of the latter. Now, \cite[Proposition 4.4]{RV21} states that the functions $\underline{\Upsilon}_p^\theta, \overline{\Upsilon}_p^\theta$ never take the values $\pm \infty$ and \cite[Proposition 4.5]{RV21} states that they are continuous as well, and are thus semi-infinite paths. Also, as stated in \cite[Proposition 4.6]{RV21}, we have $\underline{\Upsilon}_p^\theta(t)= \overline{\Upsilon}_p^\theta(t)=y$, or in other words, both the above semi-infinite paths emanate from $p$.

  To complete the proof, it now remains to show that $\inte (\underline{\Upsilon}_p^\theta)\cap \cT_\uparrow^\theta=\emptyset$ and $\inte (\overline{\Upsilon}_p^\theta)\cap \cT_\uparrow^\theta=\emptyset$, and we just show the former since the proof of the latter is analogous.  Now, one case is immediately ruled out by \cite[Lemma 5.7]{RV21} which states that if $p'=(y',t')\notin \inte(\underline{\Upsilon}_p^\theta)$, then $\Gamma_{p',\uparrow}^\theta\cap \inte(\underline{\Upsilon}_p^\theta)=\emptyset$ for any geodesic $\Gamma_{p',\uparrow}^\theta$. Further, it states that, in the remaining case when $p'\in \inte (\underline{\Upsilon}^\theta_p)$, a geodesic $\Gamma_{p',\uparrow}^\theta$ can possibly lie on $\underline{\Upsilon}^\theta_p$ for a while but does not intersect $\underline{\Upsilon}^\theta_p$ once it leaves. To complete the proof, we need only show that this case cannot occur as well, and we now do so by using precompactness properties of geodesics to reduce it to the first case which we already know how to analyse.

 Suppose that for some $p'=(y',t')\in \inte (\underline{\Upsilon}_p^\theta)$, there exists a geodesic $\Gamma^\theta_{p',\uparrow}$ such that $\Gamma^\theta_{p',\uparrow}\lvert_{[t',t'+\delta]}=\underline{\Upsilon}^\theta_p\lvert_{[t',t'+\delta]}$ for some $\delta>0$ with $t'+\delta<t$. Now, by Lemma \ref{lem:47} (2), we can find a non-empty open set $V$ such that all upward $\theta$-directed semi-infinite geodesics $\Gamma_{q,\uparrow}^\theta$ for $q\in V$ pass through $(\Gamma^\theta_{p',\uparrow},(t'+\delta/2),t'+\delta/2)\in \inte(\underline{\Upsilon}^\theta_{p})$. %
 However, this would imply that there exist points $q\notin \inte (\underline{\Upsilon}_p^\theta)$ for which $\Gamma^\theta_{q,\uparrow}$ meets $\inte(\underline{\Upsilon}^\theta_{p})$, but we already showed that this is not possible in the previous paragraph.
\end{proof}
Having shown that interfaces $\Upsilon_p^\theta$ exist for all points $p$, we now characterize the set of points admitting a unique interface.
\begin{lemma}
  \label{lem:53}
  Almost surely, for all $\theta\in \Xi_\uparrow^c$, a point $p\in \RR^2$ has a uniquely defined interface $\Upsilon^\theta_p$ if and only if $p\notin \cT^\theta_\uparrow$.
\end{lemma}
\begin{proof}
 First, we show that if $p=(y,t)$ admits two distinct interfaces $\eta_1$ and $\eta_2$, then we must have $p\in \cT^\theta_\uparrow$. Due to Lemma \ref{lem:56}, we can assume, without loss of generality, that $\eta_1$ is to the left of $\eta_2$, and we let $q\in \RR^2$ be a point in between the two. Now, since $\cI_\downarrow^\theta$ and $\cT_\uparrow^\theta$ are disjoint by definition, any geodesic $\Gamma_{q,\uparrow}^\theta$ must pass through $p$. Thus $p\in \inte (\Gamma_{q,\uparrow}^\theta)$ and this completes the proof.

 The task now is to show that any point $p\in \cT_\uparrow^\theta$ indeed admits at least two choices of the interface $\Upsilon_p^\theta$. Indeed, with the notation from Lemma \ref{lem:2}, we now show that for $p\in \cT_\uparrow^\theta$, we have $\underline{\Upsilon}_p^\theta\neq \overline{\Upsilon}_p^\theta$. Since $p\in \cT^\theta_\uparrow$, by applying Lemma \ref{lem:47} (2), we can find an open set $V$ such that for all $q\in V$ and all geodesics $\Gamma_{q,\uparrow}^\theta$, we have $p\in \inte (\Gamma^\theta_{q,\uparrow})$, As a consequence for all $q\in V$, we have $d_{p,\uparrow}^\theta(q)=0$. As a result, for any $t\in \RR$ such that $\{(x,s)\in \RR^2:s=t\}\cap V\neq \emptyset$, we have $\underline{\Upsilon}^\theta_p(t)<\overline{\Upsilon}^\theta_p(t)$ and this shows that $\underline{\Upsilon}^\theta_p$ and $\overline{\Upsilon}^\theta_p$ are distinct. This completes the proof.
\end{proof}

We now show an analogous result which in particular implies that points $p\in \cI_\downarrow^\theta$ have more than one upward $\theta$-directed geodesic emanating out from them.
\begin{figure}
  \centering
  \includegraphics[width=0.3\textwidth]{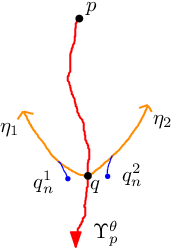}
  \caption{Proof of Lemma \ref{lem:52}-- By approximating $q$ from the left and right and taking subsequential limits of the corresponding $\theta$-directed upward semi-infinite geodesics, we obtain two distinct geodesics $\eta_1,\eta_2$ emanating from $q$.}
  \label{fig:disjoint}
\end{figure}
\begin{lemma}
  \label{lem:52}
 Almost surely, for all $\theta\in \Xi_\uparrow^c$, for any $p=(y,t)\in \RR^2$, any interface $\Upsilon^\theta_p$, any point $q=(\Upsilon^\theta_p(s),s)$ for some $s<t$, we have $q\in \NU^\theta_{0,\uparrow}$ and in fact $\underline{\Gamma}^\theta_{q,\uparrow}(s')<\Upsilon_p^\theta(s')<\overline{\Gamma}^\theta_{q,\uparrow}(s')$ for all $s'\in (s,t)$. Further, %
if $p\notin \cT^\theta_\uparrow$, we have $\underline{\Gamma}^\theta_{q,\uparrow}(t)<y<\overline{\Gamma}^\theta_{q,\uparrow}(t)$ as well.
\end{lemma}

\begin{proof}
We refer the reader to Figure \ref{fig:disjoint} for a summary of this proof. Choose sequences of real numbers $x_n^1,x_n^2,s_n$ such that $x_n^1,x_n^2\rightarrow x$ and $s_n\rightarrow s$ as $n\rightarrow \infty$ and further, $x_n^1<\Upsilon^\theta_p(s_n)<x_n^2$ for all $n$. Define the points $q_n^j=(x_n^j,s_n)$ for $j\in \{1,2\}$ and let $\Gamma^\theta_{q_n^j,\uparrow}$ be a choice of an upward $\theta$-directed semi-infinite geodesic emanating from $q_n^j$. For $j\in \{1,2\}$, we let $\eta_j$ be a subsequential limit of the geodesics $\Gamma^\theta_{q_n^j,\uparrow}$ as $n\rightarrow \infty$, which exists by Lemma \ref{prop:12}, and note that both $\eta_1,\eta_2$ are $\theta$-directed upward semi-infinite geodesics emanating from $q$. Since $\cI_\downarrow^\theta$ and $\cT_\uparrow^\theta$ are disjoint, all the geodesics $\Gamma^\theta_{q_n^1,\uparrow}$ are to the left of $\Upsilon^\theta_p$ and all the $\Gamma^\theta_{q_n^2,\uparrow}$ are to the right of $\Upsilon^\theta_p$. Thus we must have $\eta_1(s')\leq \Upsilon^\theta_p(s')\leq \eta_2(s')$ for all $s'\in (s,t)$.

  Since $\eta_1,\eta_2$ are themselves $\theta$-directed upward semi-infinite geodesics, by again using the disjointness of $\cT^\theta_\downarrow$ and $\cI^\theta_\uparrow$, we in fact obtain $\eta_1(s')<\Upsilon^\theta_p(s')<\eta_2(s')$ for all $s'\in (s,t)$. In particular, this shows that $q\in \NU^\theta_{0,\uparrow}$. %
  Further, that we have $\eta_1(t)\leq y\leq \eta_2(t)$ is clear by the above along with the continuity of $\eta_1,\eta_2$, and if there were to be equality in the above on either side, we would have $p\in \inte(\Gamma^\theta_{q,\uparrow})$ and this would imply that $p\in \cT^\theta_\uparrow$.
\end{proof}

\subsection{Non-existence of bi-infinite geodesics and the proof of Theorem \ref{thm:4}}
\label{sec:non-existence}
In this section, we provide the proof of Theorem \ref{thm:4}. As an intermediate step, we need to first rule out the presence of bi-infinite geodesics in the directed landscape. Such a statement was shown for exponential LPP in \cite{BHS22}, and we now adapt the arguments therein to obtain the corresponding statement for the directed landscape. %

\begin{figure}
  \centering
  \includegraphics[width=0.65\textwidth]{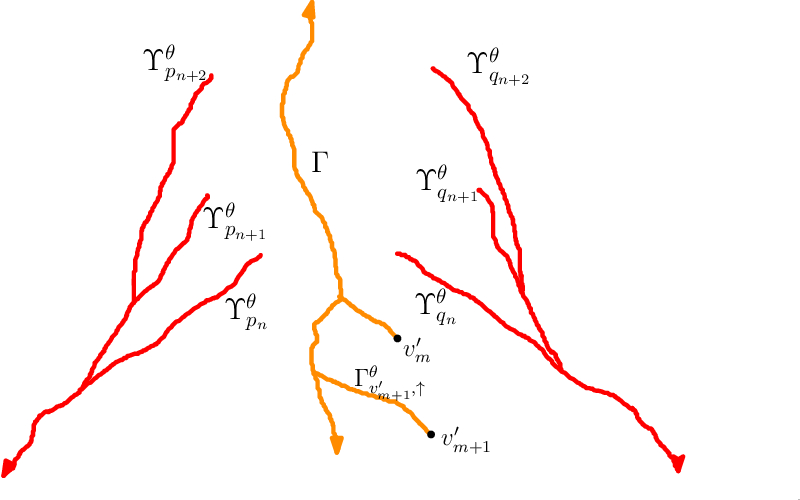}
  \caption{Proof of Theorem \ref{thm:4}: If the directed forest $\cI_\downarrow^\theta$ has multiple connected components, then we can use a compactness argument to sandwich a bi-infinite geodesic between them, but the latter does not exist.}
  \label{fig:int-tree}
\end{figure}

\begin{proposition}
  \label{prop:5}
  Almost surely, there does not exist any bi-infinite geodesic for $\cL$, that is, there does not exist any bi-infinite path $\Gamma\colon(-\infty,\infty)\rightarrow \RR$ such that $\Gamma\lvert_{[s,t]}$ is a geodesic for every $s<t$.
\end{proposition}
\begin{proof}
  We consider a bi-infinite geodesic $\Gamma$ as a concatenation of the downward and upward semi-infinite geodesics $\Gamma\lvert_{(-\infty,0]}$ and $\Gamma\lvert_{[0,\infty)}$. As discussed in Section \ref{sec:later}, due to \cite[Theorem 2.5 (i)]{BSS22}, the above paths must respectively be $\theta_\Gamma^\downarrow,\theta_\Gamma^\uparrow$ directed for some $\theta_\Gamma^\downarrow,\theta_\Gamma^\uparrow\in \RR$, and we refer to these as the downward and upward directions for the bi-infinite geodesic $\Gamma$. Now, we claim that almost surely, for all bi-infinite geodesics $\Gamma$, we must have $\theta_\Gamma^\uparrow=\theta_\Gamma^\downarrow$. Indeed, if the above were not true, then for some $c_\Gamma>0$ and all large enough $n$, we must have $|\Gamma(0)-(\Gamma(-n)+\Gamma(n))/2| \geq c_\Gamma n$ along with $\max(|\Gamma(-n)|,|\Gamma(n)|)\leq 2\max(|\theta_{\Gamma}^\downarrow|,|\theta_{\Gamma}^\uparrow|)n$. However, this is not possible since, for all $n$, the finite geodesic $\Gamma\lvert_{[-n,n]}$ must satisfy the uniform transversal fluctuation estimate from Proposition \ref{prop:trans-simul}. As a result, we must have $\theta_\Gamma^\downarrow=\theta_\Gamma^\uparrow$, and we thus define $\theta_\Gamma=\theta^\downarrow_\Gamma = \theta^\uparrow_\Gamma$ and say that $\Gamma$ is $\theta_\Gamma$-directed.

  Thus, to rule out the existence of bi-infinite geodesics and prove the proposition, it suffices to show that for any $M,K>0$, the event
  \begin{equation}
  \label{eq:58}
  \cE_{M,K}=\{\exists \textrm{ bi-infinite geodesic }\Gamma\textrm{ with } |\Gamma(0)|<K\textrm{ and }|\theta_{\Gamma}|<M\}
\end{equation}
satisfies $\PP(\cE_{M,K})=0$, and this is the goal of the remainder of the proof. To achieve this goal, we first obtain an $O_M(n^{2/3})$ bound on the expectation of the cardinality of the set $\cH_{n}^{M}$ defined by
  \begin{equation}
    \label{eq:41}
    \cH_n^{M}=\{\gamma_{(x,-n)}^{(y,n)}(0):x,y\in [-M n,M n]\},
  \end{equation}
  where $\gamma_{(x,-n)}^{(y,n)}$ is varied over all possible geodesics between $(x,-n)$ and $(y,n)$ in the above. Indeed, with $|\cH_n^{M}|$ denoting the cardinality of the set $\cH_n^{M}$, we claim that there exists a constant $C$ such that we have
  \begin{equation}
    \label{eq:49}
    \EE |\cH_n^{M}|  \leq CM^2 n^{2/3}
  \end{equation}
  for all $n\in \NN$. To see this, using the notation $[\![x,y]\!]$ to denote $[x,y]\cap \ZZ$ for $x<y$, we first define the points $a_i$ for $i\in [\![0,M n^{1/3}]\!]$ such that $a_i=-M n+2n^{2/3}i$. Having done so, we note that $\cH_n^{M}= \bigcup_{i,j\in [\![0,M n^{1/3}-1]\!]}\cH_n^{M,i,j}$, where we define
  \begin{equation}
    \label{eq:50}
    \cH_n^{M,i,j}=\{\gamma_{(x,-n)}^{(y,n)}(0):x\in [a_i,a_{i+1}],y\in [a_j,a_{j+1}]\}.
  \end{equation}
  Now, by using the skew and KPZ symmetries satisfied (Proposition \ref{prop:symm}) by the directed landscape, we obtain that
  \begin{equation}
    \label{eq:55}
    |\cH_n^{M,i,j}|\stackrel{d}{=}|\cH_n^{n^{-1/3}}|\stackrel{d}{=}|\cH_1^{1}|
  \end{equation}
  for all $i,j\in [\![0,M n^{1/3}]\!]$. Note that as a consequence of Lemma \ref{lem:45}, we know that $|\cH_1^1|$ has stretched exponential tails and thus has finite expectation. Thus, by using this along with \eqref{eq:55}, we obtain that for some constant $C$, $\EE |\cH_{n}^{M ,i,j}|\leq C$ for all $i,j\in [\![0,M n^{1/3}-1]\!]$ and this immediately yields the desired bound \eqref{eq:49}.

  We now return to the event $\cE_{M,K}$. First, we note that for $z\in \RR$ and $\alpha<\beta\in \RR$, if we use $\cA_{n,K}^{\alpha,\beta}(z)$ to denote the event that there exist $x,y\in [\alpha , \beta]$ and a geodesic $\gamma_{(x,-n )}^{(y,n )}$ satisfying $\gamma_{(x,-n )}^{(y,n )}(0)\in [z-K,z+K]$, then we must have
\begin{equation}
  \label{eq:59}
  \cE_{M,K}\subseteq \liminf_{n\rightarrow\infty}\cA_{n,K}^{-M n, M n}(0).
\end{equation}
Thus, in order to show that $\PP(\cE_{M,K})=0$, it suffices to establish that
\begin{equation}
  \label{eq:62}
  \lim_{n\rightarrow \infty}\PP(\cA_{n,K}^{-M n, M n}(0))=0,
\end{equation}
and this is the goal for the remainder of the proof. We now consider a point $z_{n}^*$ sampled, independently of $\cL$, uniformly in $[-Mn/2,Mn/2]$ and consider the event $\cA_{n,K}^{-Mn,Mn}(z_{n}^*)$. It is now easy to see that
\begin{equation}
  \label{eq:63}
  \cA_{n,K}^{-Mn,Mn}(z_{n}^*)=\{[z_{n}^*-K,z_{n}^*+K]\cap \cH_n^{M}\neq \emptyset\},
\end{equation}
and since the point $z_{n}^*$ is sampled independently of the landscape $\cL$, we immediately obtain that, for some constant $C$,
\begin{equation}
  \label{eq:54}
  \PP( \cA_{n,K}^{-Mn,Mn}(z_{n}^*))\leq \frac{2K}{M n}\EE |\cH_n^{M}|\leq CKM n^{-1/3},
\end{equation}
where we have used \eqref{eq:49} to obtain the last inequality. As a result, there must exist a deterministic $z_{n}\in [-Mn/2,Mn/2]$ for which we have
\begin{equation}
  \label{eq:64}
  \PP( \cA_{n,K}^{-Mn,Mn}(z_{n}))\leq  CKMn^{-1/3}.
\end{equation}
Now, as a consequence of the translational invariance of the directed landscape, we obtain that $\PP(\cA_{n,K }^{-Mn,Mn}(z_{n}))=\PP(\cA_{n,K}^{-M n-z_{n},M n-z_{n}}(0))$. Since $z_{n}\in [-Mn/2,Mn/2]$, we immediately obtain
\begin{equation}
  \label{eq:56}
  \cA_{n,K}^{-M n-z_{n},M n-z_{n}}(0)\supseteq \cA_{n,K}^{-Mn/2,Mn/2}(0),
\end{equation}
and as a result of this, we obtain the inequality
\begin{equation}
  \label{eq:57}
  \PP(\cA_{n,K}^{-Mn/2,Mn/2}(0))\leq CMKn^{-1/3}
\end{equation}
which in particular implies \eqref{eq:62} if we replace $M$ by $2M$. Thus we have shown that $\PP(\cE_{M,K})=0$ for all $M,K>0$, and this completes the proof.

\end{proof}
We combine the above with a compactness argument (summarised in Figure \ref{fig:int-tree}) to argue that downward interface portraits, which are a priori only forests, are actually one-ended trees simultaneously for all directions $\theta\in \Xi_\uparrow^c$. 

\begin{proof}[Proof of Theorem \ref{thm:4}]
 Suppose the contrary, and let $\theta\in \Xi_\uparrow^c$ be an angle such that there exist two points $p=(x^*,s^*),q=(y^*,t^*)$ and interfaces $\Upsilon^\theta_p,\Upsilon^\theta_q$ satisfying $\Upsilon^\theta_p(t)< \Upsilon^\theta_q(t)$ for all $t\leq \min(s^*,t^*)$. It is easy to see that we can further assume that $p,q\notin \cT_\uparrow^\theta$. %
 Now, as shown in Figure \ref{fig:int-tree}, we construct sequences of points $p_n=(x_n,s_n),q_n=(y_n,t_n)\not\in \cT_\uparrow^\theta$ for all $n\in \NN$ such that $\Upsilon^\theta_{p_n}$ is to the left of $\Upsilon^\theta_{p_m}$ for all $n>m$ and similarly $\Upsilon^\theta_{q_n}$ to the right of $\Upsilon^\theta_{q_m}$ for all $n>m$. Further, we will have that $s_n,t_n\rightarrow \infty$ as $n\rightarrow \infty$.

  For the construction, begin by defining $p_1=p,q_1=q$ and simply define $s_n=t_n=\max(s_1,t_1)+n$. Assume that $x_n,y_n$ have already been defined; we first define $y_{n+1}$, and the definition of $x_{n+1}$ will proceed analogously. Consider a geodesic $\Gamma_{q_n,\uparrow}^\theta$, and choose $y_{n+1}$ such that $y_{n+1}>\Gamma_{q_n,\uparrow}^\theta(t_{n+1})$. Now, since $\cT^\theta_\uparrow\cap \cI^\theta_\downarrow=\emptyset$, any interface $\Upsilon^\theta_{p_{n+1}}$ cannot cross $\Gamma^\theta_{q_{n},\uparrow}$ and thus must stay to the right of $\Upsilon^\theta_{q_{n}}$.
The definition of $p_{n+1}$ can be made similarly by looking at left-most instead of right-most geodesics.

Now consider a sequence of points $v'_m=(x'_m,s'_m)$ such that $s'_m\downarrow -\infty$ as $m\rightarrow \infty$ and $v'_m$ lies between the interfaces $\Upsilon^\theta_{p_1},\Upsilon^\theta_{q_1}$ and let $\Gamma^\theta_{v'_m,\uparrow}$ be a $\theta$-upward semi-infinite geodesic emanating from $v'_m$. Now given an $K\in \NN$, choose $n$ large enough so that $s_n=t_n>K$ and by using that the geodesics $\Gamma^\theta_{v_m',\uparrow}$ cannot cross the interfaces $\Upsilon^\theta_{p_n},\Upsilon^\theta_{q_n}$, we obtain by a compactness argument that there exists a subsequence $\{m_i\}$ along which $\Gamma^\theta_{v'_m,\uparrow}(-K),\Gamma^\theta_{v'_m,\uparrow}(K)$ both converge. Now, by applying Proposition \ref{prop:9}, along a possibly different subsequence $\{m_i\}$, the geodesics $\Gamma^\theta_{v'_m,\uparrow}\lvert_{[-K,K]}$ converge uniformly. By doing this for all $K\in \NN$ and taking a diagonal subsequence, we obtain that there exists a bi-infinite path $\Gamma$ such that along a possibly different subsequence $\{m_i\}$, $\Gamma^\theta_ {v'_m,\uparrow}$ converges to $\Gamma$ locally uniformly. Using that uniform limits of finite geodesics are geodesics (Proposition \ref{prop:9}), we obtain that $\Gamma\lvert_{[s,t]}$ for any $s<t$ is a geodesic and thus $\Gamma$ is a bi-infinite geodesic. However, this cannot happen by Proposition \ref{prop:5} and thus the assumption of there existing $\Upsilon^\theta_p,\Upsilon^\theta_q$ which are disjoint must be false. This completes the proof. %
\end{proof}

We note that an argument of the above flavour, where one takes subsequential limits of certain geodesics to obtain a bi-infinite geodesic, has recently appeared \cite{AHH22} in the area of first passage percolation.

\section{The duality in the continuum and its consequences}
\label{sec:dual}
\subsection{The prelimiting landscapes $\cL^n$}
For this section, we define the vectors $\bv=(1,1)$ and $\bw=(1,-1)$. We use $X^n={X^n_{(i,j)}}$ to denote a family of i.i.d.\ $\exp(1)$ random variables indexed by $(i,j)\in \ZZ^2$, and these variables serve as the weights of the prelimiting exponential LPP model which will converge to the directed landscape. For any two points $p=(x_1,y_1),q=(x_2,y_2)\in \ZZ^2$ with $p\leq q$ in the sense that $x_1\leq x_2$ and $y_1\leq y_2$, we define the passage time $T^n(p,q)$ by
\begin{equation}
  \label{eq:30}
  T^n(p,q)=\max_{\pi:p\rightarrow q}\ell(\pi),
\end{equation}
where the maximum is over all up-right lattice paths $\pi$ from $p$ to $q$ with the corresponding weight $\ell(\pi)$ defined by $\ell(\pi)\coloneqq\sum_{v\in \pi\setminus \{p\}}X_v$\footnote{Contrary to a convention often used in the literature, we do not add the weight $X^n_p$ of the first vertex $p$ here. In practice, with regards to results, there is no difference between the two conventions and in particular, all the results in this section are true for both conventions. The convenient aspect of not including the first vertex in the length is that $T^n(p,q), T^n(q,r)$ are independent for every $p<q<r$ and satisfy an exact composition law similar to \eqref{eq:9}.}. The almost surely unique lattice path attaining the maximum in \eqref{eq:30} is called a geodesic and is denoted by $\gamma_{p,q}^{n,\dis}$. We also think of $\gamma_{p,q}^{n,\dis}$ as a function such that $m\mathbf{v}+\gamma_{p,q}^{n,\dis}(m)\mathbf{w}\in \gamma_{p,q}^{n,\dis}$, where the above makes sense for $m$ lying in a finite subset of $(1/2)\ZZ$. To introduce the prelimiting directed landscapes $\cL^n$ which are the objects which will actually converge to the landscape $\cL$, we need some additional notation. For a given $p\in \ZZ^2$, we define the set $\boxx(p)$ by
\begin{equation}
  \label{eq:8}
  \boxx(p)=p+\{s\mathbf{v}+x\mathbf{w}: s\in (-1/4,1/4], x\in (-1/2,1/2]\},
\end{equation}
and note that $\RR^2$ is the disjoint union $\bigsqcup_{p\in \ZZ^2}\boxx(p)$. For $q\in \RR^2$, we use the notation $\mathfrak{r}(q)$ to denote the unique $p\in \ZZ^2$ for which $q\in \boxx(p)$.
We now define the prelimiting directed landscapes $\cL^n$ such that for any $(x,s),(y,t)$ with $s<t$,
\begin{equation}
  \label{eq:31}
  \cL^n(x,s;y,t)=2^{-4/3}n^{-1/3}\left(T^n(\mathfrak{r}( sn\bv+2^{2/3}xn^{2/3}\bw),\mathfrak{r}( tn\bv+2^{2/3}yn^{2/3}\bw))-4(t-s)n\right).
\end{equation}
For $p=(x_1,s_1),q=(x_2,s_2)\in \RR^2$ with $s_1<s_2$, the corresponding rescaled geodesic $\gamma_{p,q}^{n,\res}\colon (s_1,s_2)\rightarrow \RR$ is defined by
\begin{equation}
  \label{eq:1}
  \gamma_{p,q}^{n,\res}(s)=2^{-2/3}n^{-2/3}\gamma^{n,\dis}_{\mathfrak{r}(2^{2/3}x_1n^{2/3}\mathbf{w}+s_1n \mathbf{v}), \mathfrak{r}(2^{2/3}x_2n^{2/3}\mathbf{w}+s_2n \mathbf{v})}(ns)
\end{equation}
for $ns\in (1/2)\ZZ$ and interpolating linearly for intermediate values of $s$. The following result from \cite{DV21} %
 shows that $\cL^n$ converges to $\cL$ as $n\rightarrow \infty$ and that the corresponding finite geodesics converge as well.
\begin{proposition}[{\cite[Theorem 1.7, Theorem 1.8, Remark 1.10]{DV21}}]
  \label{prop:6}
  There is a coupling between $X^n$ for all $n\in \NN$ such that $\cL^n\rightarrow \cL$ almost surely with respect to the locally uniform topology on $\RR^4_\uparrow$. Further, in the above coupling, almost surely, for any points $p_n\rightarrow p\in \RR^2$ and $q_n\rightarrow q\in \RR^2$, the geodesics $\gamma_{p_n,q_n}^{n,\res}$ are precompact in the uniform topology and any subsequential limit is a geodesic $\gamma_p^q$.
\end{proposition}
Just as we have semi-infinite geodesics in the directed landscape, they also exist in the prelimit \cite{FMP09,FP05}. For $p\in \ZZ^2$, we will use the notation $\Gamma^{n,\dis}_{p,\downarrow},\Gamma^{n,\dis}_{p,\uparrow}$ to denote such geodesics for the LPP given by the weights $X^n$ and going in the direction $-\bv,\bv$ respectively, and these are almost surely unique for all $p$. Thus, the semi-infinite geodesic $\Gamma_{p,\uparrow}^{n,\dis}\subseteq \ZZ^2$ (resp.\ $\Gamma_{p,\downarrow}^{n,\dis}\subseteq \ZZ^2$) is an up-right (resp.\ down-left) semi-infinite lattice-path emanating from $p$ and going asymptotically in the direction $\bv$ (resp.\ $-\bv$) in the sense that if we have a sequence $(x_m,y_m)$ consisting of pairwise distinct elements such that $(x_m,y_m)\in \Gamma_{p,\uparrow}^{n,\dis}$ (resp.\ $\Gamma_{p,\uparrow}^{n,\dis}$) for all $m$, then a.s.\ we have $\lim_{m\rightarrow \infty}x_m=\infty$ (resp.\ $-\infty$) and $\lim_{m\rightarrow \infty} x_m/y_m=1$ (resp.\ $1$). Apart from being viewed as a subset of $\ZZ^2$, such a geodesic will also be thought of as a function such that the point 
\begin{equation}
  \label{eq:44}
 m\bv+\Gamma_{p,\uparrow}^{n,\dis}(m)\bw 
\end{equation}
 lies on the geodesic, where the above makes sense for all large enough $m\in(1/2)\ZZ$, with a similar statement being true for $\Gamma^{n,\dis}_{p,\downarrow}$. We note that the quantities $\Gamma^{n,\dis}_{p,\uparrow}(m)/m$ and $\Gamma^{n,\dis}_{p,\downarrow}(m)/m$ converge to $0$ as $m$ goes to $\infty$ and $-\infty$ respectively since $\Gamma^{n,\dis}_{p,\downarrow},\Gamma^{n,\dis}_{p,\uparrow}$ go asymptotically in the $\bv,-\bv$ direction respectively. As per our usual convention, we will simply use $\Gamma_{p}^{n,\dis}$ to denote $\Gamma_{p,\downarrow}^{n,\dis}$.

For $p=(x,s)\in \RR^2$, we define the rescaled geodesic $\Gamma^{n,\res}_{p}=\Gamma^{n,\res}_{p,\downarrow}$ by
\begin{equation}
  \label{eq:34}
  \Gamma^{n,\res}_{p,\downarrow}(t)=2^{-2/3}n^{-2/3}\Gamma^{n,\dis}_{\mathfrak{r}(2^{2/3}xn^{2/3}\bw+sn\bv),\downarrow}(nt).
\end{equation}
for $nt\in (1/2)\ZZ$ and interpolate linearly in between.
A similar definition allows us to define $\Gamma^{n,\res}_{p,\uparrow}$. %
As in the continuum, the geodesics $\{\Gamma^{n,\dis}_{p}\}_{p\in \ZZ^2}$ exhibit coalescence \cite{FP05} and form a discrete one-ended tree which we call $\cT_\downarrow(T^n)$. Also, we use the notation $\cT_\downarrow(\cL^n)$ to denote the tree obtained by considering the graphs of the scaled geodesics \eqref{eq:34} as $p$ varies. %

 Due to the above coalescence, we can define the prelimiting Busemann functions $\cB_\downarrow^{n,\dis}(p,q)$ for $p,q\in \ZZ^2$ by
\begin{equation}
  \label{eq:32}
  \cB_\downarrow^{n,\dis}(p,q)=T^n(z,p)-T^n(z,q),
\end{equation}
where $z\in \ZZ^2$ denotes the point where the geodesics $\Gamma^{n,\dis}_p$ and $\Gamma^{n,\dis}_q$ first meet. Similarly, for $p,q\in \RR^2$, we can define
\begin{equation}
  \label{eq:35}
  \cB_\downarrow^{n,\res}(p,q)=\cL^n(z;p)-\cL^n(z;q),
\end{equation}
where $z$ is the first meeting point of $\Gamma_p^{n,\res}$ and $\Gamma_q^{n,\res}$. We note that the upward objects $\cB_\uparrow^{n,\dis},\cB_\uparrow^{n,\res}$ can be defined similarly.

Though the coupling in Proposition \ref{prop:6}, a priori, only guarantees the convergence of finite geodesics, it is in fact true that both prelimiting infinite geodesics and prelimiting Busemann functions converge to their continuum counterparts, and we state such a result now.
\begin{proposition}
  \label{prop:1}
  The coupling of $\cL^n$ with $\cL$ from Proposition \ref{prop:6} can be chosen such that the following precompactness for discrete semi-infinite geodesics holds. Almost surely, for any sequence of points $p_n\rightarrow p\in \RR^2$, the geodesics $\Gamma_{p_n,\downarrow}^{n,\res}$ (resp.\ $\Gamma_{p_n,\uparrow}^{n,\res}$) are precompact in the locally uniform topology and any subsequential limit is a geodesic $\Gamma_{p,\downarrow}$ (resp.\ $\Gamma_{p,\uparrow}$). Further, almost surely, the Busemann functions $\cB_\downarrow^{n,\res}(\cdot,\0)$ (resp.\ $\cB_\uparrow^{n,\res}(\cdot,\0)$) converge locally uniformly to their continuum counterpart $\cB_\downarrow(\cdot,\0)$ (resp.\ $\cB_\uparrow(\cdot,\0)$).
\end{proposition}
From this point in the paper, until Lemma \ref{lem:18}, we will consider the $X^n$ to be coupled in the above manner. The proof of Proposition \ref{prop:1} involves a technical and tedious coupling argument which we defer to the appendix. We now come to discrete interface portraits and their rescaled variants. For exponential LPP, one can define the interface portrait  $\cI_\downarrow(T^n)$ by the requirement that a given dual edge $e^*\in (\ZZ^2)^*=\ZZ^2+(1/2,1/2)$ belongs to $\cI_\downarrow(T^n)$ if the primal edge $e$ crossing the dual edge does not belong to the geodesic tree $\cT_\uparrow(T^n)$. The interface portrait $\cI_\downarrow(T^n)$, which is a priori only a forest, was shown to a.s.\ be a tree in \cite{Pim16}, and for each point $p\in (\ZZ^2)^*$, one can define the interface $\Upsilon^{n,\dis}_p=\Upsilon^{n,\dis}_{p,\downarrow}$ as the a.s.\ unique down-left semi-infinite lattice path emanating from $p$ in the tree $\cI_\downarrow(T^n)$. Just as we think of the geodesics $\Gamma_p^{n,\dis}$ as functions (see \eqref{eq:44}), we similarly think of discrete interfaces as functions such that $m\bv+\Upsilon^{n,\dis}_{p}(m)\bw$ lies on the interface $\Upsilon_p^{n,\dis}$, where the above makes sense for all negative enough $m\in(1/2)\ZZ+(1/4)$.

One can similarly define upward interfaces and interface portraits and we denote these by $\Upsilon^{n,\dis}_{p,\uparrow}$ and $\cI_\uparrow(T^n)$. Since we will be interested in the continuum limit, we define rescaled prelimiting interfaces by defining for $p\in (x,s)\in \RR^2$,
\begin{equation}
  \label{eq:341}
  \Upsilon^{n,\res}_{p,\downarrow}(t)=2^{-2/3}n^{-2/3}\Upsilon^{n,\dis}_{\mathfrak{r} (2^{2/3}xn^{2/3}\bw+sn\bv),\downarrow}(nt),
\end{equation}
for $nt\in (1/2)\ZZ+(1/4)$, and we extend this to arbitrary $t$ by linearly interpolating in between. We denote the corresponding tree, formed by the graphs of the $\Upsilon_{p,\downarrow}^{n,\res}$ for different values of $p$, by $\cI_\downarrow(\cL^n)$. As always, we have a corresponding $\uparrow$ version of the above objects, and in case we are working with the $\downarrow$ version, we omit $\downarrow$ from the notation. Before moving to the next section, we state a useful result regarding the convergence of discrete interfaces to continuum interfaces.

\begin{figure}
  \centering
  \captionsetup[subfigure]{labelformat=empty}
  \setbox9=\hbox{\includegraphics[width=0.4\textwidth]{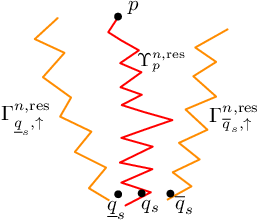}}
  \subcaptionbox{}{\raisebox{\dimexpr\ht9-\height}{\includegraphics[width=0.4\textwidth]{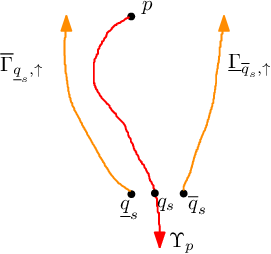}}}
  \hfill
  \subcaptionbox{}{\includegraphics[width=0.4\textwidth]{intface-conv1}}
  \caption{Proof of Lemma \ref{lem:23}: First, for points $\underline{q}_s, \overline{q}_s$ within $\varepsilon$ distance to the left and right of $q_s$, we construct geodesics $\overline{\Gamma}_{\underline{q}_s,\uparrow}$ and $\underline{\Gamma}_{\overline{q}_s,\uparrow}$ which pass strictly to the left and right of $p$ respectively. Due to the convergence of discrete rescaled semi-infinite geodesics to continuum ones, the same must hold true for all large $n$ for the geodesics $\Gamma_{\underline{q}_s,\uparrow}^{n,\res}$ and $\Gamma_{\overline{q}_s,\uparrow}^{n,\res}$ and as a result, the discrete interface $\Upsilon_p^{n,\res}$, which is sandwiched between the above two geodesics, must be $2\varepsilon$-close to $q_s$ at time $s$.}
  \label{fig:intface-conv}
\end{figure}

\begin{lemma}
  \label{lem:23}
  Almost surely, for every $p=(y,t)\notin \cT_\uparrow$, the discrete interfaces $\Upsilon^{n,\res}_{p}$ converge to $\Upsilon_p$ pointwise  as $n\rightarrow \infty$.
\end{lemma}

\begin{proof}
We refer the reader to Figure \ref{fig:intface-conv} for a summary of the proof. Recall that, since $p\notin \cT_\uparrow$, there is a unique interface $\Upsilon_p$ as a consequence of Lemma \ref{lem:53}. Fix $s<t$ and define $q_s=(\Upsilon_p(s),s)$. Consider the left-most and right-most geodesics $\underline{\Gamma}_{q_s,\uparrow},\overline{\Gamma}_{q_s,\uparrow}$ and note that by Lemma \ref{lem:52}, we have
  \begin{equation}
    \label{eq:39}
   \underline{\Gamma}_{q_s,\uparrow}(t)<y<\overline{\Gamma}_{q_s,\uparrow}(t).
  \end{equation}
  For a given $\varepsilon>0$, consider the points $\underline{q}_s=(\Upsilon_p(s)-\varepsilon,s),\overline{q}_s=(\Upsilon_p(s)+\varepsilon, s)$ and note that by semi-infinite geodesic ordering (Lemma \ref{prop:2}) along with \eqref{eq:39}, we have
  \begin{equation}
    \label{eq:47}
    \overline{\Gamma}_{\underline{q}_s,\uparrow}(t)<y<\underline{\Gamma}_{\overline{q}_s,\uparrow}(t).
  \end{equation}
  Since any subsequential limits of the prelimiting geodesics $\Gamma^{n,\res}_{\underline{q}_s,\uparrow},\Gamma^{n,\res}_{\overline{q}_s,\uparrow}$ are (Proposition \ref{prop:1}) $0$-directed upward semi-infinite geodesics from $\underline{q}_s,\overline{q}_s$ respectively, we obtain by \eqref{eq:47} that there exists a $\delta>0$ such that for all large enough $n$, we have
  \begin{equation}
    \label{eq:48}
    \Gamma^{n,\res}_{\underline{q}_s,\uparrow}(t)<y-\delta<y+\delta<\Gamma^{n,\res}_{\overline{q}_s,\uparrow}(t).
  \end{equation}
Thus, since discrete interfaces and semi-infinite geodesics cannot cross each other, the above implies that $|\Upsilon^{n,\res}_{p}(s)-\Upsilon_p(s)|<2\varepsilon$ for all large enough $n$. Since $\varepsilon$ was arbitrary, this shows that $\Upsilon_p^{n,\res}(s)\rightarrow \Upsilon_p(s)$ as $n\rightarrow \infty$, thereby completing the proof.
 \end{proof}

\subsection{Constructing the duality coupling}
To obtain the duality for the directed landscape, we will in fact take a subsequential limit of the duality \cite{Pim16} in the discrete prelimiting model of exponential LPP along with Proposition \ref{prop:6}, the recently established convergence of exponential LPP to the directed landscape.  To begin, we recall the prelimiting duality from \cite{Pim16}. %
\begin{proposition}[{\cite[Lemma 2]{Pim16}}]
  \label{prop:7}
Given the i.i.d.\ $\exp(1)$ field $X^n=\{X^n_{(i,j)}\}_{(i,j)\in \ZZ^2}$, consider the dual field $\wX^n=\{\wX^n_{(i^*,j^*)}\}_{(i^*,j^*)\in (\ZZ^2)^*}$ defined by
  \begin{equation}
    \label{eq:33}
    \wX^n_{(i+1/2,j+1/2)}=\min \left(\cB_\downarrow^{n,\dis}( (i,j+1),\0),\cB_\downarrow^{n,\dis}( (i+1,j),\0)\right)-\cB_\downarrow^{n,\dis}( (i,j),\0),
  \end{equation}
  where we identify the dual $(\ZZ^2)^*$ with $\ZZ^2+(1/2,1/2)$. Using $\wT^n$ to denote passage times for the LPP model defined by $\wX^n$, consider the geodesic tree $\cT_\uparrow(\wT^n)\subseteq (\ZZ^2)^*$ and the interface portrait $\cI_\downarrow(\wT^n)\subseteq ((\ZZ^2)^*)^*=\ZZ^2$. Then we have the a.s.\ equalities $\cT_\downarrow(T^n)=\cI_\downarrow(\wT^n)$, $\cI_\uparrow(T^n)=\cT_\uparrow(\wT^n)$. %
\end{proposition}
We note that in \cite[Lemma 2]{Pim16}, the above result is instead simply stated as $\cI_\downarrow(\wT^n)\stackrel{d}{=} \cT_\downarrow(T^n)$, but the proof proceeds precisely by constructing $\widetilde{T}^n$ via \eqref{eq:33}. Indeed, if we use $\widetilde \cB_\uparrow^{n,\dis}$ to denote the upward Busemann function with respect to the weights $\widetilde X^n$, then the latter weights are defined exactly so that the function $z^*\mapsto \cB_\uparrow^{n,\dis}(z^*,(1/2,1/2))$ defined for $z^*\in (\ZZ^2)^*$ is equal to $-B^{\downarrow*}$, where the object $B^{\downarrow*}$ appears in the proof of \cite[Lemma 2]{Pim16}. As an example of a setting where the duality from \cite{Pim16} is phrased using the weights $\widetilde X^n$ similarly to Proposition \ref{prop:7}, we refer the reader to \cite[Lemma 4.2, Lemma 4.3]{FMP09}.

We note that as in Proposition \ref{prop:7}, we will work with the case of $0$-directed geodesic trees and interface portraits throughout this section. Later, when we come to the proofs of Theorems \ref{thm:20} and \ref{thm:6}, we shall use the skew-invariance of the directed landscape (Proposition \ref{prop:symm}) to transfer results from the $0$-directed case to the case of general directions $\theta\in \RR$.

Now, just as we defined prelimiting directed landscapes $\cL^n$ using $T^n$ in \eqref{eq:31}, we can analogously define the prelimiting directed landscapes $\widetilde{\cL}^n$ such that for any $(x,s),(y,t)$ with $s<t$, we have
\begin{equation}
  \label{eq:31d}
  \widetilde{\cL}^n(x,s;y,t)=2^{-4/3}n^{-1/3}\left(\widetilde{T}^n(\mathfrak{r}( sn\bv+2^{2/3}xn^{2/3}\bw),\mathfrak{r}( tn\bv+2^{2/3}yn^{2/3}\bw)-4(t-s)n\right).
\end{equation}
 As an immediate consequence of Proposition \ref{prop:7}, we now have the following lemma.
\begin{lemma}
  \label{lem:19}
  Almost surely, we have the exact equalities $\cT_{\downarrow}(\cL^n)=\cI_{\downarrow}(\wcL^n)$ and $\cI_{\uparrow}(\cL^n)=\cT_{\uparrow}(\wcL^n)$ for all $n$.
\end{lemma}
Our aim is to show the analogue of the above lemma for the directed landscape. We begin by noting the following important consequence of Propositions \ref{prop:6}, \ref{prop:1}.
\begin{lemma}
  \label{lem:18}
  The sequence $(\cL^n,\wcL^n)$ is tight in $n$ for the locally uniform topology on $\RR_\uparrow^4$, and thus there exists a subsequential weak limit $(\cL,\wcL)$ where $\wcL\stackrel{d}{=}\cL$. Further, we can couple $X^n,\wX^n,\cL,\wcL$ such that, along a deterministic subsequence $\cN\subseteq \NN$, $\cL^n\rightarrow \cL$ and $\wcL^n\rightarrow\wcL$ almost surely locally uniformly as $n\rightarrow \infty$ with the additional discrete semi-infinite geodesic precompactness from Proposition \ref{prop:1} holding for both $\cL^n$ and $\widetilde \cL^n$.
\end{lemma}

For the remainder of this section, we will work with the coupling from the above lemma, and will often have convergence results which hold only along the subsequence $\cN$. To keep the notation clean, if a sequence $a_n$ converges to a limit $a$ along the subsequence $\cN$, we will simply say that $a_n$ converges along $\cN$ to $a$.

For the coupling $(\cL,\wcL)$ to have any utility, we need to obtain a continuum analogue of Lemma \ref{lem:19}, and we shall obtain this by tracking how the equality in Lemma \ref{lem:19} behaves as $n\rightarrow \infty$. We now show that in a certain long-term sense, the prelimiting geodesic tree $\cT_\downarrow(\cL^n)$ converges along $\cN$ to $\cT_\downarrow(\cL)$.

To introduce the notion of convergence that we will use, we return to the general setting of directed forests introduced at the end of Section \ref{sec:thm2pf}. Given a random downward directed forest $\fF$, we say that $\cT_\downarrow(\cL^n)\rightarrow \fF$ a.s.\ in the long-term sense if for any $a<b\in \RR$ and any bounded open set $U\subseteq \RR^2$ with $\overline{U}\subseteq \{(x,t):t>b\}$, $\bigcup_{(x,t)\in U}\Gamma^n_{(x,t)}\cap (\RR\times [a,b])$ a.s.\ converges to $ \bigcup_{(x,t)\in U}\fF_{(x,t)}\cap (\RR\times [a,b])$ in the Hausdorff metric as $n\rightarrow \infty$. A similar definition can be made for the convergence of prelimiting interface portraits to downward/upward directed forests as well. We now prove a scaling limit result for discrete downward $0$-directed geodesic trees.
\begin{lemma}
   \label{lem:21}
   Almost surely, along the subsequence $\cN$ and in the long-sense, $\cT_{\downarrow}(\cL^n)$ converges to $\cT_{\downarrow}(\cL)$ and $\cT_{\uparrow}(\widetilde{\cL}^n)$ converges to $\cT_{\uparrow}(\widetilde{\cL})$.
 \end{lemma}

 \begin{proof}
 We only show the a.s.\ long-term sense convergence of $\cT_{\downarrow}(\cL^n)$ to $\cT_{\downarrow}(\cL)$, and the proof of the other convergence is analogous.  Let $a<b\in \RR$ and the bounded open set $U$ be as in the definition of long-term convergence. By using basic transversal fluctuation estimates for semi-infinite geodesics (Lemma \ref{lem:3}), we can first find rational points $p,q$ and geodesics $\Gamma_p,\Gamma_q$ such that $U$ lies in between $\Gamma_p$ and $\Gamma_q$. By using the ordering of semi-infinite geodesics (Lemma \ref{prop:2}), we obtain that all geodesics $\Gamma_{(y,t)}$ for $(y,t)\in U$ must in fact lie between $\Gamma_p,\Gamma_q$. Thus, by applying Lemma \ref{lem:45}, we obtain that there exist finitely many geodesics $\gamma_i\colon (a,b)\rightarrow \RR$ with $1\leq i\leq i_0$ such that for any $(y,t)\in U$, we have $\Gamma_{(y,t)}\lvert_{(a,b)}=\gamma_i$ for some $1\leq i \leq i_0$.

   Now suppose that for some $\delta>0$, and for infinitely many $n\in \cN$, there exist points $(y_n,t_n)\in U$ such that $\Gamma^{n,\res}_{(y_n,t_n)}\lvert_{(a,b)}\not \subseteq B_\delta(\cup_{1\leq i\leq i_0}\gamma_i\lvert_{(a,b)})$, where the latter denotes a spatial $\delta$-neighbourhood of the set $\cup_{1\leq i\leq i_0}\gamma_i\lvert_{(a,b)}$. By Proposition \ref{prop:1}, this would imply on taking a limit that there exists a geodesic $\Gamma_{(y,t)}$ such that $\Gamma_{(y,t)}\lvert_ {(a,b)}\neq \gamma_i$ for all $1\leq i \leq i_0$ and this would be a contradiction. Thus we have shown that for any fixed $\delta>0$, almost surely for all large enough $n\in \cN$, we have
   \begin{equation}
     \label{eq:20}
     \bigcup_{(y,t)\in U} \Gamma^{n,\res}_{(y,t)}\lvert_{(a,b)}\subseteq B_\delta(\bigcup_{1\leq i\leq i_0}\gamma_i\lvert_{(a,b)}).
   \end{equation}

  To complete the proof, we now need to establish that for each fixed $\delta>0$, almost surely, for each $i$ with $1\leq i\leq i_0$, there do exist infinitely many $n\in \cN$ and points $(y_n,t_n)\in  U$ for which we have $\Gamma^{n,\res}_{(y_n,t_n)}\lvert_{(a,b)}\subseteq B_\delta(\gamma_i\lvert_{(a,b)})$. To do so, we first recall that by definition, we can find a point $(y,t)\in U$ for which we have $\Gamma_{(y,t)}\lvert_{(a,b)}=\gamma_i$. Now, since $U$ is open, the set $\inte (\Gamma_{(y,t)})\cap U$ is non-empty, and we choose a point $(y',t')$ lying in this set. Now, if there exists more than one geodesic $\Gamma_{(y',t')}$, then we would have a geodesic bubble, and thus by Theorem \ref{thm:2}, the geodesic $\Gamma_{(y',t')}$ must be a.s.\ unique. Also, we note that $\Gamma_{(y',t')}\lvert_{(a,b)}=\gamma_i$. With the above at hand, we simply define $(y_n,t_n)\in U$ such that $(y_n,t_n)\rightarrow (y',t')$ as $n\rightarrow \infty$. Finally, by an application of Proposition \ref{prop:1} and the a.s.\ uniqueness of the geodesic $\Gamma_{(y',t')}$, we obtain that the geodesics $\Gamma_{(y_n,t_n)}^{n,\res}$ converge locally uniformly along $\cN$ to the geodesic $\Gamma_{(y',t')}$, and this immediately implies the statement mentioned at the beginning of the paragraph. This completes the proof. %
 \end{proof}

With the above convergence statements in place, we can now complete the proof of Theorem \ref{thm:20}.
\begin{proof}[Proof of Theorem \ref{thm:20}]
  We first take $\theta=0$ and show that $(\cL,\wcL)$ is the required coupling in the sense that we have the almost sure equalities $\cI_\uparrow(\cL)= \cT_\uparrow(\wcL)$ and $\cT_\downarrow(\cL)=\cI_\downarrow(\wcL)$. We only show the former equality and the latter follows by analogous arguments.

  We first show that $\cI_\uparrow(\cL)\subseteq \cT_\uparrow(\wcL)$ almost surely. If this were not the case, then there would exist an interface $\Upsilon_{p,\uparrow}$ such that $\inte (\Upsilon_{p,\uparrow})\not \subseteq \cT_\uparrow(\wcL)$.
Now, for any point $q\in \inte(\Upsilon_{p,\uparrow})$, we know that $q\notin \cT_\downarrow(\cL)$ since $\cI_\uparrow(\cL)$ and $\cT_\downarrow(\cL)$ are disjoint by definition. Thus by Lemma \ref{lem:53}, there is a unique interface $\Upsilon_{q,\uparrow}$ and considering the discrete rescaled interfaces $\Upsilon^{n,\res}_{q,\uparrow}$ introduced in \eqref{eq:341} and by using %
  Lemma \ref{lem:23}, we obtain that a.s.\ $\Upsilon^{n,\res}_{q,\uparrow}\rightarrow \Upsilon_{q,\uparrow}$ pointwise along $\cN$. However, we also have $\Upsilon_{q,\uparrow}^{n,\res}\subseteq \cI_{\uparrow}(\cL^n)=\cT_\uparrow(\wcL^n)$ which, by Lemma \ref{lem:21}, converges along $\cN$ to $\cT_\uparrow(\wcL)$ in the long term sense, thereby implying that $\inte (\Upsilon_{q,\uparrow})\subseteq \cT_\uparrow(\wcL)$. Since $q$ was an arbitrary point satisfying $q\in \inte(\Upsilon_{p,\uparrow})$, this shows that $\inte(\Upsilon_{p,\uparrow})\subseteq \cT_\uparrow(\widetilde \cL)$, contrary to what we had assumed earlier.

  It now remains to show that $\cT_\uparrow(\wcL)\subseteq \cI_\uparrow(\cL)$ almost surely, and before this, we  quickly show that $\cT_\uparrow(\wcL)$ and $\cT_\downarrow(\cL)$ are a.s.\ disjoint. Indeed, if the above disjointness were not true, then we would have a cycle in $\cT_\uparrow(\wcL)$ by an application of Lemma \ref{lem:53} and the fact that $\cI_\uparrow(\cL)\subseteq \cT_\uparrow(\wcL)$, but this is not possible since, by Corollary \ref{thm:3}, we already know that $\cT_\uparrow(\wcL)$ is a.s.\ a one-ended tree.

  Now, with the goal of eventually obtaining a contradiction, assume that $\cT_\uparrow(\wcL)\not \subseteq \cI_\uparrow(\cL)$. Then there exists a $0$-directed upward semi-infinite geodesic $\widetilde{\Gamma}_{p,\uparrow}$ for $\wcL$ with $\inte (\widetilde{\Gamma}_{p,\uparrow})\not\subseteq\cI_\uparrow(\cL)$. As a consequence of Theorem \ref{thm:2}, for any $q=(x,s)\in \widetilde{\Gamma}_{p,\uparrow}$ with $q\neq p$, $\widetilde{\Gamma}_{p,\uparrow}\lvert_{[s,\infty)}$ is the unique $0$-directed upward semi-infinite geodesic $\widetilde{\Gamma}_{q,\uparrow}$ emanating from $q$. We can further choose $q$ so that $\widetilde{\Gamma}_{q,\uparrow}\not \subseteq \cI_\uparrow(\cL)$, and we also know that $q\notin \cT_\downarrow(\cL)$ due to the disjointness proved in the previous paragraph. However, due to the uniqueness of $\widetilde{\Gamma}_{q,\uparrow}$ along with the discrete semi-infinite geodesic precompactness from Lemma \ref{lem:18} and Proposition \ref{prop:1}, we know that a.s.\ $\widetilde{\Gamma}_{q,\uparrow}^{n,\res}\rightarrow \widetilde{\Gamma}_{q,\uparrow}$ pointwise along $\cN$. Since, by Lemma \ref{lem:19}, $\widetilde{\Gamma}^{n,\res}_{q,\uparrow}=\Upsilon_{q,\uparrow}^{n,\res}$ and because $q\notin\cT_\downarrow(\cL)$, we have by Lemma \ref{lem:23} that $\Upsilon_{q,\uparrow}^{n,\res}$ must converge a.s.\ along $\cN$ to $\Upsilon_{q,\uparrow}$ pointwise. This shows that $\widetilde{\Gamma}_{q,\uparrow}=\Upsilon_{q,\uparrow}$, and since $q$ was an arbitrary point on $\inte(\widetilde{\Gamma}_{p,\uparrow})$, this contradicts the assumption that $\inte (\widetilde{\Gamma}_{p,\uparrow})\not \subseteq \cI_\uparrow(\cL)$. This completes the proof of Theorem \ref{thm:20} for the case $\theta=0$.

  Now, to construct the coupling $(\cL_\theta,\wcL_\theta)$ for any fixed $\theta\in \RR$, we note that by Lemma \ref{lem:skew}, on applying an appropriate $\theta$-dependent skew transformation to both the landscapes $\cL,\wcL$, we obtain a coupled pair of landscapes $(\cL_\theta,\wcL_\theta)$ such that the $\theta$-directed geodesic trees and interface portraits corresponding to the latter are simple skew transformations of the $0$-directed geodesic trees and interface portraits corresponding to $(\cL,\wcL)$. This yields the desired coupling $(\cL_\theta,\wcL_\theta)$ and completes the proof.

\end{proof}

\subsection{Using the duality to investigate fractal aspects of the geodesic tree}
The aim of this section is to prove Theorem \ref{thm:6}. We begin with using the H\"older continuity properties of geodesics to obtain their Hausdorff dimension as a subset of $\RR^2$.
\begin{lemma}
  \label{lem:4}
  For any fixed $p\in \RR^2$, almost surely, every segment of $\Gamma_p$ has Hausdorff dimension $4/3$ as a subset of $\RR^2$.
\end{lemma}
\begin{proof}
  Due to the translational invariance (Proposition \ref{prop:symm}) of the directed landscape, we can take $p=\0$ without loss of generality. We begin with showing that $\dim \Gamma_\0\leq 4/3$ almost surely. To see this, we first note that with the geodesic web $\cW=\bigcup_{(x,s;y,t)\in \QQ_\uparrow^4}\inte(\gamma_{(x,s)}^{(y,t)})$ defined in Proposition \ref{lem:50}, it suffices to show that $\dim \cW\leq 4/3$ almost surely. By using the stability of Hausdorff dimension with respect to countable unions, it equivalently suffices to show that for any fixed $(x,s;y,t)\in \RR_\uparrow^4$, we have $\dim \gamma_{(x,s)}^{(y,t)}\leq 4/3$ almost surely. Now, by Proposition \ref{prop:3}, the geodesic $\gamma\coloneqq \gamma_{(x,s)}^{(y,t)}$ is a.s.\ $2/3-\varepsilon$ H\"older continuous for any fixed $\varepsilon>0$, and by a standard argument (see e.g.\ \cite[Lemma 5.1.2]{BP16}), this  yields the desired $2-2/3=4/3$ upper bound on the Hausdorff dimension.

We now come to the lower bound. Since every segment of $\Gamma_\0$ contains a segment $\Gamma_\0\lvert_{[s_0,t_0]}$ where $s_0,t_0$ are rational, it suffices to show that for any \emph{fixed} $s_0<t_0\leq 0$, $\dim \Gamma_\0\lvert_{[s_0,t_0]}\geq 4/3$ almost surely, and we equivalently show that for any fixed $\varepsilon>0$, we a.s.\ have $\dim \Gamma_\0\lvert_{[s_0,t_0]} \geq 4/3-\varepsilon$. To do so, we first consider the occupation measure $\zeta$ of $\Gamma_\0$ defined such that for any measurable $A\subseteq \RR^2$, we have
  \begin{equation}
    \label{eq:26}
    \zeta(A)=\leb( t\leq 0: (\Gamma_0(t),t)\in A),
  \end{equation}
  where $\leb$ denotes the Lebesgue measure. In other words, $\zeta$ is the pushforward of the Lebesgue measure on $(-\infty,0]$ under the map $t\mapsto (\Gamma_\0(t), t)$. We note that $\zeta$ is a measure supported on the graph of the geodesic $\Gamma_\0$ and since $\Gamma_\0$ is a semi-infinite path, we have $\zeta( [-M,M]\times [s_0,t_0])>0$ almost surely for all $M$ large enough such that $\Gamma_\0\lvert_{[s_0,t_0]}\cap ([-M/2,M/2]\times [s_0,t_0])\neq \emptyset$. In fact, due to the above and the mass distribution principle \cite[Lemma 1.2.8]{BP16}, in order to show that $\dim \Gamma_\0\lvert_{[s_0,t_0]} \geq 4/3-\varepsilon$, it suffices to show that almost surely, for any fixed $M\in \NN$,
  \begin{equation}
    \label{eq:27}
    \sup_{A\subseteq [-M,M]\times [s_0,t_0]} \frac{\zeta(A)}{\diam(A)^{4/3-\varepsilon}}<\infty,
  \end{equation}
  where the supremum above is over all balls $A$ in $[-M,M]\times [s_0,t_0]$ and $\diam (A)$ refers to the diameter of $A$ with respect to the usual Euclidean metric on $\RR^2$. %

  We now show \eqref{eq:27} and the main tool that we use for this is Proposition \ref{prop:4}. By standard arguments, it suffices to show \eqref{eq:27} with $A$ being restricted to be dyadic squares instead of balls. To be precise, using the notation $[\![a,b]\!]$ to denote the set $[a,b]\cap \ZZ$ for $a<b\in \RR$, we need only show that for some a.s.\ finite positive random variable $R$ and for all $k\in \NN$ and $i\in [\![-M2^k,M2^k-1]\!]$, $j\in [\![s_02^k,t_02^k-1]\!]$, we have
  \begin{equation}
    \label{eq:29}
    \zeta([i2^{-k},(i+1)2^{-k}]\times [j2^{-k},(j+1)2^{-k}])\leq R 2^{-(4/3-\varepsilon)k}.
  \end{equation}
  To show the above, we use Proposition \ref{prop:4} along with a simple Borel-Cantelli argument. Indeed, Proposition \ref{prop:4} and a union bound immediately imply that for some constants $C,c$ and all $k\in \NN$, we have
  \begin{equation}
    \label{eq:38}
    \PP( \sup_{i,j}\zeta([i2^{-k},(i+1)2^{-k}]\times [j2^{-k},(j+1)2^{-k}])>2^{-(4/3-\varepsilon)k}) \leq C\exp(-c2^{\varepsilon k}) \times (2^{k+1}M)\times 2^k(t_0-s_0).
  \end{equation}
  We now observe that the right hand side above is summable in $k$ and thus by applying the Borel-Cantelli lemma, we immediately obtain \eqref{eq:29}. This completes the proof.
\end{proof}

Now, by applying \ref{lem:8}, we use the above to obtain the dimension of the geodesic tree.
\begin{lemma}
  \label{lem:31}
  For any fixed direction $\theta\in \RR$, the geodesic tree $\cT^\theta_\downarrow$ a.s.\ has Hausdorff dimension $4/3$.
\end{lemma}
\begin{proof}
  As in the proof of Theorem \ref{thm:20}, by an application of Lemma \ref{lem:skew}, it suffices to prove the statement for $\theta=0$. By Lemma \ref{lem:4}, almost surely simultaneously for all rational $p\in \RR^2$, we know that $\inte (\Gamma_p)$ has Hausdorff dimension $4/3$, and the result now follows by Lemma \ref{lem:8}. %
\end{proof}

Finally, we use that there are countably many confluence points (Lemma \ref{lem:44}) to obtain that the geodesic tree has countably many trifurcation points, where by a trifurcation point of a directed forest $\fF$, we mean a point where two semi-infinite paths $\fF_p,\fF_q$ starting at some $p\neq q\in \RR^2$ first meet.
\begin{lemma}
  \label{lem:46}
 For any fixed $\theta\in \RR$, the tree $\cT^\theta_\downarrow$ a.s.\ has only countably many trifurcation points.
\end{lemma}

\begin{proof}
 Again, by Lemma \ref{lem:skew}, it suffices to specialise to the case $\theta=0$. Since each trifurcation point is a confluence point in the sense of Lemma \ref{lem:44}, we obtain that there are a.s.\ at most countably many trifurcation points. To show that there are at least infinitely many such points, we consider the geodesics $\{\Gamma_{(n,0)}\}_{n\in \NN}$ and note that any point where two such semi-infinite geodesics first meet is a trifurcation point of $\cT_\downarrow$. We already showed at the end of the proof of Lemma \ref{lem:44} that there are infinitely many such first meeting points, and the proof is thus complete.
\end{proof}

We can now complete the proof of Theorem \ref{thm:6}.
\begin{proof}[Proof of Theorem \ref{thm:6}]
  By Theorem \ref{thm:20}, we know that $\cI_\uparrow^\theta(\cL_\theta)=\cT_\uparrow^\theta(\wcL_\theta)$ holds almost surely, and by Lemma \ref{lem:31} and Lemma \ref{lem:46}, we know that $\dim\cT_\uparrow^\theta(\wcL_\theta)=4/3$ and that $\mathrm{Tri}(\cT_\uparrow^\theta(\wcL_\theta))$ is countable almost surely. Thus it remains to show that $\NU_0^\theta(\cL_\theta)=\cI^\theta_\uparrow(\cL_\theta)$ almost surely, and that the set of points $p$ admitting three distinct geodesics $\Gamma^\theta_p$ is equal to $\mathrm{Tri}(\cI_\uparrow^\theta(\cL_\theta))$. We first show the former.

  Note that by Lemma \ref{lem:52}, every point on the interior of any interface $\Upsilon^\theta_{p,\uparrow}$ admits two $\theta$-directed downward semi-infinite geodesics and this shows that $\cI_\uparrow^\theta(\cL_\theta)\subseteq \NU_0^\theta(\cL_\theta)$, and we now need to show the reverse inclusion. For this, let $p=(x,s)\in \NU_0^\theta(\cL_\theta)=\NU_1^\theta(\cL_\theta)$ and consider the left-most and right-most geodesics $\underline{\Gamma}^\theta_p,\overline{\Gamma}^\theta_p$. Let $q=(y,t)\in \QQ^2$ with $t<s$ be a point satisfying $\underline{\Gamma}_p^\theta(t)<y<\overline{\Gamma}_p^\theta(t)$. Consider the interface $\Upsilon^\theta_{q,\uparrow}$, which is unique since $q\in \QQ^2$ and $\QQ^2\cap \cT_\downarrow^\theta(\cL_\theta)=\emptyset$ a.s.\ by Lemma \ref{lem:1}. Note that since $\cI^\theta_\uparrow(\cL_\theta)$ and $\cT^\theta_\downarrow(\cL_\theta)$ are disjoint, the interface $\Upsilon^\theta_{q,\uparrow}$ must pass through $p$, and this shows that $p\in \cI_\uparrow^\theta(\cL_\theta)$, thereby showing that a.s.\ $\NU^\theta_0(\cL_\theta)\subseteq \cI_\uparrow^\theta(\cL_\theta)$ and completing the proof of the equality.

  We now show the corresponding statement for the set of trifurcation points $\mathrm{Tri}(\cI_\uparrow^\theta(\cL_\theta))$. Locally, we define the set $\three$ by
  \begin{equation}
    \label{eq:19}
    \three = \{ p\in \RR^2: \textnormal{there exist three distinct geodesics } \Gamma^\theta_p\},
  \end{equation}
and our goal now is to show that $\three= \mathrm{Tri}(\cI_\uparrow^\theta(\cL_\theta))$ almost surely. To obtain $\three\subseteq \mathrm{Tri}(\cI_\uparrow^\theta(\cL_\theta))$, we use a similar argument as above. Let $p=(x,s)$ be a point in the former set and let $\eta_1,\eta_2,\eta_3$ denote three choices of the geodesic $\Gamma^\theta_p$, ordered from left to right. Note that $\eta_1(t)\neq \eta_2(t)\neq \eta_3(t)$ for all $t\in (s-\varepsilon,s)$ for some random $\varepsilon>0$ as a consequence of Theorem \ref{thm:2} and the fact that all $\eta_1,\eta_2,\eta_3$ eventually coalesce (Proposition \ref{prop:11}). Thus we can find rational points $q_1=(y_1,t)$ and $q_2=(y_2,t)$ such that $\eta_1(t)<y_1<\eta_2(t)<y_2<\eta_3(t)$ for some $t\in (-\infty,s)$. Now we consider the interfaces $\Upsilon^\theta_{q_1,\uparrow},\Upsilon^\theta_{q_2,\uparrow}$ and use that $\cT^\theta_\downarrow(\cL_\theta)$ and $\cI^\theta_\uparrow(\cL_\theta)$ are disjoint to obtain that both the above interfaces must pass through $p$, and that $p$ must be the first point at which they intersect. This implies that $p\in \mathrm{Tri}(\cI_\uparrow^\theta(\cL_\theta))$, thereby establishing the inclusion $\three\subseteq \mathrm{Tri}(\cI_\uparrow^\theta(\cL_\theta))$.

For the reverse inclusion, consider a point $p=(y,t)\in \mathrm{Tri}(\cI_\uparrow^\theta(\cL_\theta))$ and let $\Upsilon^\theta_{p_1,\uparrow},\Upsilon^\theta_{p_2,\uparrow}$ be two choices of interfaces, intersecting for the first time at $p$, with the former being to the left of the latter. Now choose a sequence of rational points $y^1_n,y^2_n,y^3_n,t_n$ such that $t_n$ increases to $t$ and $y^1_n<\Upsilon^\theta_{p_1,\uparrow}(t_n)<y^2_n<\Upsilon^\theta_{p_2,\uparrow}(t_n)<y^3_n$. Define the points $q^j_n=(y_n^j,t_n)$ for $j\in \{1,2,3\}$ and let $\eta_j$ be a subsequential limit of the geodesics $\Gamma^\theta_{q^j_n}$ as $n\rightarrow \infty$, and note that by Lemma \ref{prop:12}, $\eta_1,\eta_2,\eta_3$ must all be $\theta$-directed downward semi-infinite geodesics emanating from $p$. If we write $p_1=(x_1,s_1)$ and $p_2=(x_2,s_2)$, then note that we must have $\eta_1(s')< \Upsilon_{p_1,\uparrow}^\theta(s')< \eta_2(s')< \Upsilon_{p_2,\uparrow}^\theta(s')< \eta_3(s')$ for all $s'\in (\max\{s_1,s_2\},t)$ since $\cT^\theta_\downarrow(\cL_\theta)$ and $\cI^\theta_\uparrow(\cL_\theta)$ are disjoint. Thus $\eta_1,\eta_2,\eta_3$ are distinct $\theta$-directed downward semi-infinite geodesics emanating from $p$, and we have $p\in \three$, thereby showing the inclusion $\mathrm{Tri}(\cI_\uparrow^\theta(\cL_\theta))\subseteq \three$ and completing the proof. 
\end{proof}

\section{Appendix}
\label{sec:appendix}
In this short appendix, we provide the proof of Proposition \ref{prop:1}. We note that the arguments here are heavily based on the corresponding arguments in \cite{DV21}. Regarding notation, for a function $f\colon \RR\rightarrow \RR$, and $s<t$, we will use $\cL^n(f,s;y,t)$ and $\cL(f,s;y,t)$ to denote $\sup_{x\in \RR}(f(x)+\cL^n(x,s;y,t))$ and $\sup_{x\in \RR}(f(x)+\cL(x,s;y,t))$ respectively. We now come to the proof of the proposition.
\begin{proof}[Proof of Proposition \ref{prop:1}]
We only exhibit a coupling which works for downward infinite geodesics and downward Busemann functions; the argument for obtaining convergence for the corresponding upward objects is analogous. For this proof, we use $B^{n,k}$ to denote the function defined by $B^{n,k}(x)= \cB^{n,\res}_\downarrow((x,k),(0,k))$ and define $B^k$ by $B^k(x)=\cB_\downarrow( (x,k), (0,k) )$. Consider the sequence $(\cL^n, \{B^{n,k}\}_{k\in \ZZ})$. We already know by Proposition \ref{prop:6} that $\cL^n\stackrel{d}{\rightarrow}\cL$ with respect to the local uniform topology. Further, it is a fact (see e.g.\ \cite[Theorem 4.2]{Sep17}) that for each $k\in \ZZ$, and for any fixed $p_0\in \{x+y=2k\}$, the discrete process $x\mapsto \cB^{n,\dis}_\downarrow ( (x, 2k-x), p_0)$ for $x\in \ZZ$ is a random walk with its jump distribution being the difference of two independent $\exp(1/2)$ variables. As a consequence, the rescaled Busemann function $B^{n,k}(x)$ is a random walk rescaled by Brownian scaling and thus converges as $n\rightarrow \infty$, in distribution with respect to the locally uniform topology, to a Brownian motion of diffusivity $2$.

  Thus, since each marginal converges in distribution as $n\rightarrow \infty$, the sequence $(\cL^n, \{B^{n,k}\}_{k\in \ZZ})$ is tight, and therefore admits a subsequential limit, which we denote as $(\cL, \{\ff_k\}_{k\in \ZZ})$. For simplicity, we take a coupling so that the above convergence occurs almost surely. The goal now is to identify the law of the above limit as that would imply that the entire sequence converges, as opposed to just converging along a subsequence.

  It suffices to identify the joint law of $\{\ff_{-K},\cdots,\ff_{K}\}$ and $\cL$ for any fixed $K\in \NN$. We know that each $\ff_i$ is marginally a Brownian motion of diffusivity $2$ and further, we know that each $\ff_i$ is independent of $\cL\lvert_{ \{(x,s;y,t)\in \RR_\uparrow^4, s>i\}}$, since the above is true before taking the limit $n\rightarrow \infty$ as well.
  Now note that by the metric composition law, for all $n$, we have the identity
  \begin{equation}
    \label{eq:4}
    B^{n,k}(x)=\cL^n(B^{n,-K},-K;x,k)-\cL^n(B^{n,-K},-K;0,k).
  \end{equation}
  Recall that we already know that the pair $(B^{n,-K},\cL^n)$ converges a.s.\ locally uniformly along a subsequence to $(\ff_{-K},\cL)$. By a simple Borel Cantelli argument, we also know that for some random $a_1,a_2>0$ and all large enough (random) $n$, we have $B^{n,-K}(x)\leq a_1|x|+a_2$ for all $x\in \RR$. Now, we invoke \cite[Theorem 6.5]{RV21} which states that under the above conditions, we can take the limit $n\rightarrow \infty$ in \eqref{eq:4} to obtain that
\begin{equation}
    \label{eq:5}
    \ff_k(x)=\cL(\ff_{-K},-K;x,k)-\cL(\ff_{-K},-K;0,k)
  \end{equation}
  holds a.s.\ for all $x\in \RR$. %

   Now, in view of \eqref{eq:5}, we note that the joint law of $\{\ff_{-K},\cdots, \ff_{K}\}$ and $\cL$ described above is in fact the same as the joint law of $\{B^{-K},\dots ,B^{K}\}$ and $\cL$. Indeed, $B^{-K}$ is a Brownian motion of diffusivity $2$ (see e.g.\ \cite[Corollary 3.22]{RV21}) independent of $\cL\lvert_{ \{(x,s;y,t)\in \RR_\uparrow^4, s>-K\}}$ and furthermore,
  \begin{equation}
    \label{eq:6}
    B^k(x)= \cL(B^{-K},-K;x,k)-\cL(B^{-K},-K;0,k)
  \end{equation}
  holds due to the composition law of the directed landscape. This identifies the joint law of $\{\ff_{-K},\cdots,\ff_{K}\}$ and $\cL$. Thus, we can choose the coupling of the $X^n$ such that apart from the conditions in Proposition \ref{prop:6}, we additionally have that $B^{n,k}$ a.s.\ converges locally uniformly to $B^k$ as $n\rightarrow \infty$ for each $k\in \ZZ$.

 We now show that in the above coupling, the Busemann functions $\cB^{n,\res}_\downarrow(\cdot,\0)$ converge to $\cB_\downarrow(\cdot,\0)$ locally uniformly almost surely. Given a compact set $A\subseteq \RR^2$, we first fix $K\in \NN$ such that the line $\{t=-K\}$ lies below $A$ and note that it suffices to prove that $\cB_\downarrow^{n,\res}(\cdot, (0,-K) )$ converges uniformly a.s.\ to $\cB_\downarrow(\cdot, (0,-K) )$ on $A$. We note that for $(x,s)\in A$,
  \begin{equation}
    \label{eq:11}
    \cB^{n,\res}_\downarrow( (x,s), (0,-K) )= \cL^n( B^{n,-K},-K;x,s),
  \end{equation}
  and we already know that $B^{n,-K}$ and $\cL^n$ both converge locally uniformly a.s.\ to their respective limits. Finally, as in \eqref{eq:5}, by using the convergence statement \cite[Theorem 6.5]{RV21},
 we conclude the a.s.\ uniform convergence of \eqref{eq:11} for $(x,s)\in A$ to
  \begin{equation}
    \label{eq:14}
    \cB_\downarrow( (x,s), (0,-K) )= \cL( B^{-K},-K;x,s),
  \end{equation}
  and this finishes the proof of the convergence of Busemann functions.
  
  The final task is to obtain the convergence of downward semi-infinite geodesics; the argument for this is along the exact same lines as the corresponding argument for finite geodesics \cite[Theorem 8.5]{DV21}, and we now quickly go through it without furnishing all the details. Let $p_n=(y_n,t_n)\rightarrow p=(y,t)\in \RR^2$ be a sequence of points. Fix $K\in \NN$ large enough such that $-K<t_n$ for all $n$, and let $A$ be a compact set which is large enough such that $(\Gamma_{p}(s),s)\in \inte (A)$ for all choices of the above geodesic and all $s\in (-K,t]$. To argue that such an $A$ exists, we simply note that by using semi-infinite geodesic ordering (Lemma \ref{prop:2}) along with Lemma \ref{lem:3}, we can find two geodesics $\Gamma_{q_1},\Gamma_{q_2}$ starting from rational points $q_1,q_2$ such that $p$ lies in between $\Gamma_{q_1},\Gamma_{q_2}$.

  We now consider the sequence of sets $\Gamma^{n,\res}_{p_n}\cap A$, which we note is precompact in the Hausdorff topology since the Hausdorff topology on compact subsets of a compact space is compact. Now, by using the a.s.\ locally uniform convergence of $(\cB^{n,\res}_\downarrow(\cdot,\0),\cL^n)$ to $(\cB_\downarrow(\cdot,\0),\cL)$, we note that any subsequential limit $\cS$ of the above sets contains $p$ and is further also a geodesic set (c.f.\ \cite[Definition 5.4]{DV21}) in the sense that for any $(x,s),(y,t)\in \cS$ with $s<t$, we have
  \begin{equation}
    \label{eq:17}
    \cB_\downarrow( (y,t),\0)=\cB_\downarrow( (x,s),\0)+ \cL(x,s;y,t).
  \end{equation}
  It can be shown (see \cite[Proposition 5.5]{DV21}) that any geodesic set as above is a subset of a downward $0$-directed semi-infinite geodesic $\Gamma_{p}$. However, the set $A$ was chosen specifically to satisfy $(\Gamma_{p}(s),s)\in \inte (A)$ for $s\in (-K,t]$, and thus $\cS\subseteq \inte (A)$ as well. Thus, by the connectedness of $\Gamma^{n,\res}_{p_n}$, we obtain that $\Gamma^{n,\res}_{p_n}\lvert_{(-K,t_n]} \subseteq \inte (A)$ for all large $n$, and thus in fact $\Gamma^{n,\res}_{p_n}\lvert_{(-K,t_n]}$ a.s.\ subsequentially converges in the Hausdorff metric to a connected geodesic set $\cS'\subseteq \inte (A)$. It is not difficult to see (see \cite[Proposition 8.3]{DV21}) that a connected geodesic set $\cS'$ as in the above must be equal to $\Gamma_{p}\lvert_{(-K,t]}$ for some geodesic $\Gamma_{p}$ and thus we obtain that $\Gamma^{n,\res}_{p_n}\lvert_{(-K,t_n]}$ converges uniformly along a subsequence to $\Gamma_{p}\lvert_{(-K,t]}$. Finally, by sending $K\rightarrow \infty$ and employing a straightforward diagonal argument, we obtain a geodesic $\Gamma_p$ such that the sequence $\Gamma^{n,\res}_{p_n}$ converges subsequentially locally uniformly to $\Gamma_p$ as $n\rightarrow \infty$. This completes the proof.

\end{proof}

\printbibliography
\end{document}